\documentclass{article}
\usepackage[utf8]{inputenc}
\usepackage{amsmath,amssymb,amsthm}
\usepackage{hyperref}
\usepackage{enumitem}   
\usepackage{setspace}
\usepackage{multirow}
\usepackage{outlines}
\usepackage[ margin=1in]{geometry}
\usepackage[affil-it]{authblk}

\usepackage{graphicx}

\usepackage{float}
\usepackage{caption}
\usepackage{subcaption}

\usepackage[numbers]{natbib}
\let\cite=\citep

\begin{document}

\title{Parameter estimation of a two state delay differential equation modeling the human respiratory system}

\author{Nirjal Sapkota%
  \thanks{Electronic address: \texttt{nxs167030@utdallas.edu}; Corresponding author}}
\affil{Department of Mathematical Sciences,\\ The University of Texas at Dallas\\ Richardson, TX, 75080, USA}

\author{Janos Turi%
  \thanks{Electronic address: \texttt{turi@utdallas.edu}}}
\affil{Department of Mathematical Sciences,\\ The University of Texas at Dallas\\ Richardson, TX, 75080, USA}

\date{}

	\maketitle

	\section*{Abstract}

	We study parameter estimation for the two state model which describes the balance equation for carbon dioxide and oxygen in human respiratory system. These are nonlinear parameter dependent and because of the transport delay in the respiratory control system, they are modeled with delay differential equation. Numerically simulated noisy data are  generated and several examples are studied with  Levenberg–Marquardt and Trust-region algorithms to determine the values of unknown parameters.

	\section{Introduction}
	\label{ch:intro}

We have the human respiratory system examined in \cite{cooke1994stability,kollar2005numerical,khoo1982factors,batzel2000stability,sapkota2022stability} as

\begin{equation} \label{eq:sys1}
	\begin{aligned}
		\frac{\mathrm{d} x}{\mathrm{d} t} &= 1 - \alpha \,  V(x(t-\tau), y(t-\tau))\, x(t)\\
		\frac{\mathrm{d} y}{\mathrm{d} t} &= 1 - \beta \, V(x(t-\tau), y(t-\tau))\, y(t)\\	
	\end{aligned}
\end{equation}
where the ventilation function is given by 
\begin{equation*} \label{eq:ventp0}
	V(x(t-\tau), y(t-\tau)) = 0.14\,e^{-0.05(100-y(t-\tau))} \, x(t-\tau)
\end{equation*}

It is observed  that the stability of the equilibrium of this model depends on the parameters $\alpha$ and $\beta.$ We would like to get the estimates of these parameters from the measured (probably noisy) data. Hartung and Turi \cite{hartung2013parameter} studied the parameter identification of a two dimensional model representing the partial pressure of the respiratory control system.  While they wrote the code for the numerical scheme  to approximate the solutions of the delay differential equations and minimization of  the objective function based on Trust-region technique, we will use two algorithms and the built in functions available in Matlab \cite{MATLAB:2020}.

Parameter identification of delay differential equations are difficult theoretical problems (see \cite{verduyn2001parameter, zhang2006parameter, hartung2000parameter, nakagiri1995unique}). The question on whether there is a unique model that fits the measurement data is a important one which needs to be studied.
\section{Parameter Estimation}

Let's suppose that our parameters $\alpha$ and $\beta$ are not known. We will not consider $\tau$ to be a parameter in this study. So our goal is to find the estimates of $\alpha$ and $\beta,$ when we have the observed data specified at certain times $t_1, t_2,...t_M$. The most common approach is to minimize the least squares criterion for fitting a model to data. 

Let $x(t_i;\alpha,\beta)$ and $y(t_i;\alpha,\beta)$ be the model prediction of system (\ref{eq:sys1}) and $X_i$ and $Y_i$ be the observed data for $x$ and $y$ at the time $t_i$ for $i = 1,...,M$ observations, then the objection function is given by

\begin{equation} \label{eq:cost1}
	J(\alpha, \beta) = \sum_{i=1}^{M}\left(x(t_i;\alpha,\beta) - X_i\right)^2 + \sum_{i=1}^{M}\left(y(t_i;\alpha,\beta) - Y_i\right)^2	
\end{equation}
 
 Hartung et al \cite{hartung2000parameter, hartung1997differentiability} and Rihan \cite{rihan2021delay} have studied parameter identification and  convergence properties of numerical schemes of approximate solutions of parameter estimation problems.

\subsection{Computation of Estimates}
There are many algorithms using the iterative technique for minimizing the nonlinear objective function $J(\alpha, \beta).$ We will use two commonly used procedures called Levenberg–Marquardt algorithm \cite{levenberg1944method, marquardt1963algorithm, more1978levenberg} and Trust-region algorithm \cite{coleman1996interior, coleman1994convergence}.

These algorithms are discussed in detail in book  \cite{martins2021engineering} where for the Trust-region algorithm, the authors used the procedure from Nocedal and Wright \cite{nocedal1999numerical} with the parameters recommended from Conn et al \cite{conn2000trust}. For detailed information on optimization techniques we refer the reader to \cite{boyd2004convex, gill2019practical, kochenderfer2019algorithms, antoniou2007practical, griva2009linear}.

To find the global best-fit parameters, we should choose a starting point of the parameters close enough to the global minimum. A good inital staring point also speeds up the minimization procedure.

Next, we show some numerical examples to demonstrate this process for the human respiratory system (\ref{eq:sys1}).  The experimental data with some noise are numerically simulated. We compare the two algorithms in how they perform.

\section*{Example 1}
In this example, we generate measurements of the system (\ref{eq:sys1}) with a set of parameter values of  $ \alpha = 0.5 $, $ \beta = 0.8 $ and $ \tau = 1.$ The measurements are taken over the interval $ [T_0, T] = [0, 5]. $ Then random measurement noise which has a normal distribution with a mean of zero and a standard deviation of 0.20 are added. 
\begin{equation}
	\begin{aligned}
	x_i &= x(t_i; \alpha = 0.5, \beta = 0.8) + \mathcal{N}(0, 0.20)\\	
	y_i &= y(t_i; \alpha = 0.5, \beta = 0.8) + \mathcal{N}(0, 0.20)	
	\end{aligned}
\end{equation}

   We consider the parameters $ \alpha $ and $ \beta $ to be unknown, and the goal is to estimates these parameter values using the measurements. For this example, we start with the initial starting point of the parameters at $	\alpha = 0.3 $ and $ \beta = 0.5.	$

\subsection*{Computaion of Estimates with Levenberg–Marquardt algorithm}

\begin{figure}[H]
	\centering
	\begin{subfigure}[b]{0.48\textwidth}
		\centering
		\includegraphics[width=1\textwidth]{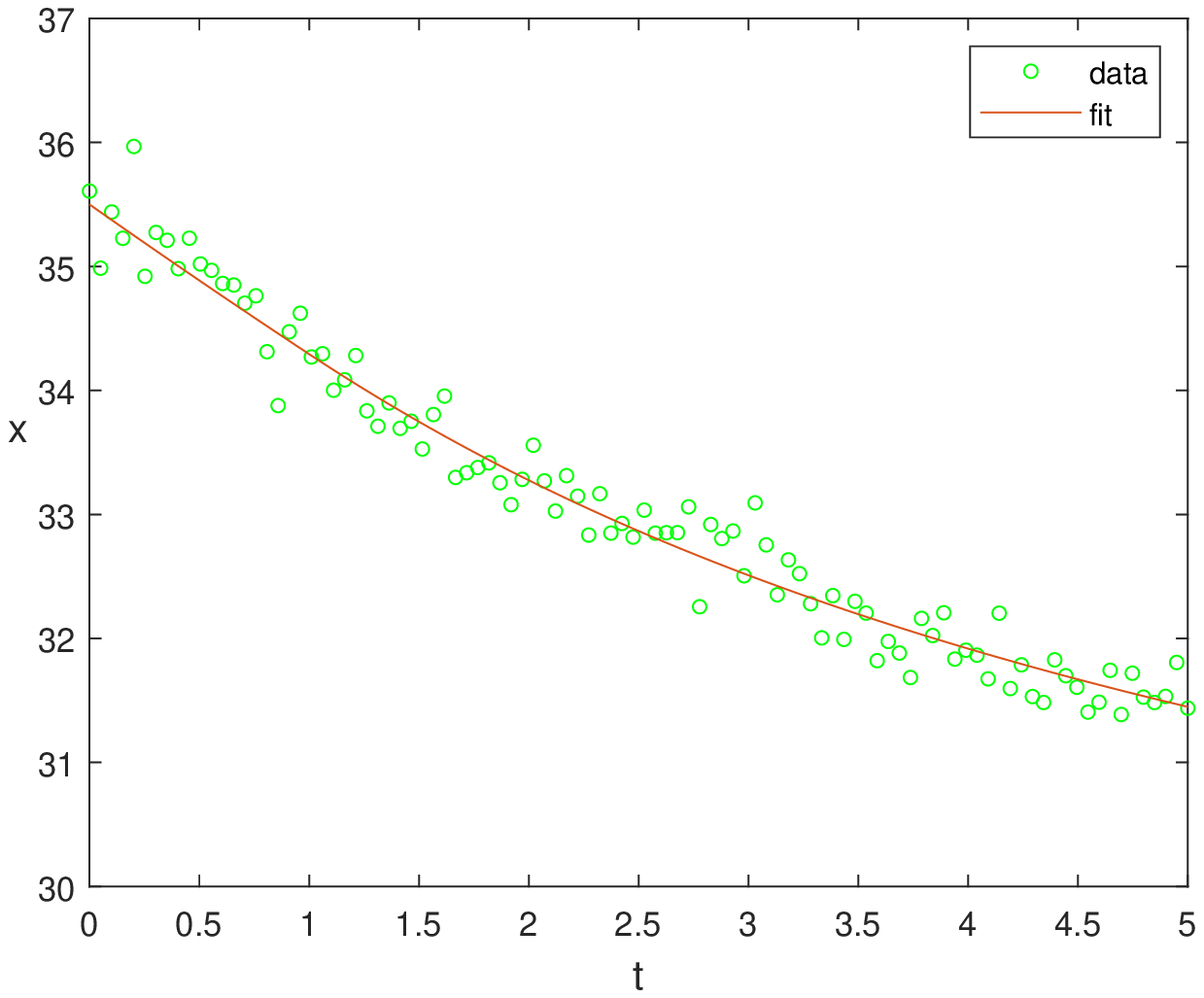}		
		\label{fig:para_est1_lmx_alpha0p3beta0p5}
	\end{subfigure}
	\hfill
	\begin{subfigure}[b]{0.48\textwidth}
		\centering
		\includegraphics[width=1\textwidth]{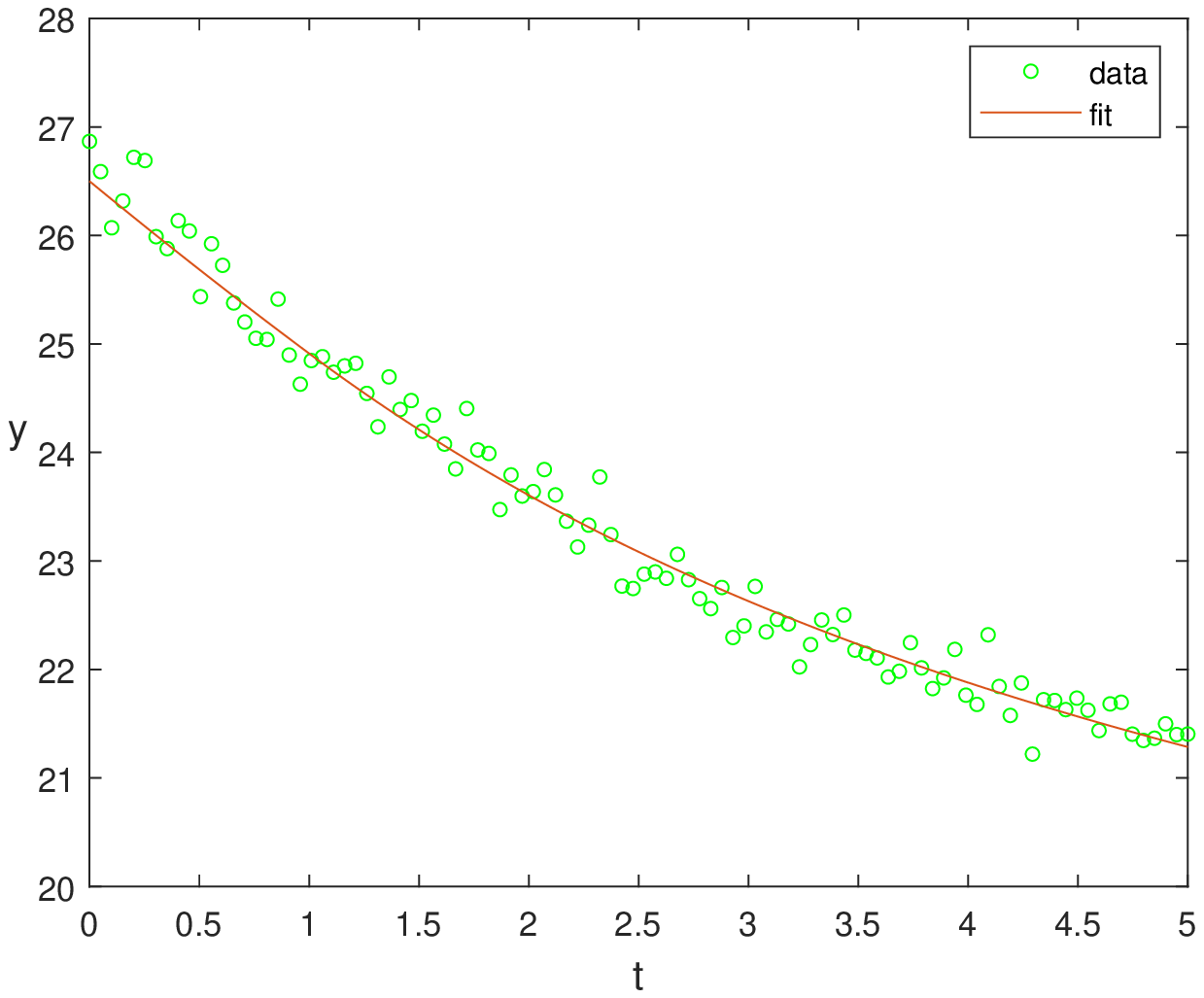}		
		\label{fig:para_est1_lmy_alpha0p3beta0p5}
	\end{subfigure}
	\caption{Data and curve fit for $x$ and $y$ with Levenberg–Marquardt algorithm for Example 1}
	\label{fig:para_est1_lm_alpha0p3beta0p5}
\end{figure}

\begin{table}[H]
	\begin{center}
		\caption{Estimation of $\alpha$ and $\beta$ with  Levenberg–Marquardt algorithm for Example 1}
		\label{table:ex1_lm}
		\begin{tabular}{cccccr}
			\hline
			Iteration &  \multicolumn{1}{p{2cm}}{\centering Function\\ count} & Residual &  \multicolumn{1}{p{2cm}}{\centering First-order\\ optimality} & Lambda & Norm of step    \\
			\hline
			0         & 3         & 738.83   & 2.02e+03   & 0.01   &             \\
			1         & 6         & 14.9085  & 139        & 0.001  & 0.326315    \\
			2         & 9         & 9.39073  & 0.977      & 0.0001 & 0.0348424   \\
			3         & 12        & 9.39033  & 0.000298   & 1e-05  & 0.000312417 \\
			4         & 15        & 9.39033  & 7.96e-07   & 1e-06  & 3.28038e-08 \\			
			\hline
		\end{tabular}
	\end{center}
\end{table}

The data points, the curve fit are plotted in Figure \ref{fig:para_est1_lm_alpha0p3beta0p5}.  Table \ref{table:ex1_lm} contains the values of  iteration count, function count, residual, first-order optimality, lambda and norm of the step. Function count is the number of function evaluations.
Lambda is the Lagrange multiplier. We used the default settings of   the Levenberg–Marquardt algorithm in Matlab. The optimization stopped because the relative norm of the current step is less than step tolerance of 1.000000e-06. After four iterations it gives the best fit of the parameters.

\begin{table}[H]
	\begin{center}
		\caption{True and best fit of  $\alpha$ and $\beta$ for Example 1}
		\label{table:ex1_true_fit_lm}
		\begin{tabular}{|c|c|c|c|}
			\hline
			&Initial & True & Best fit \\
			 \hline
		$\alpha$ & 0.3 & 0.5  & 0.5021  \\ \hline
		$\beta$  & 0.5 & 0.8 & 0.7996 \\
			\hline
		\end{tabular}
	\end{center}
\end{table}

\begin{figure}[H]
	\centering
	\begin{subfigure}[b]{0.48\textwidth}
		\centering
		\includegraphics[width=1\textwidth]{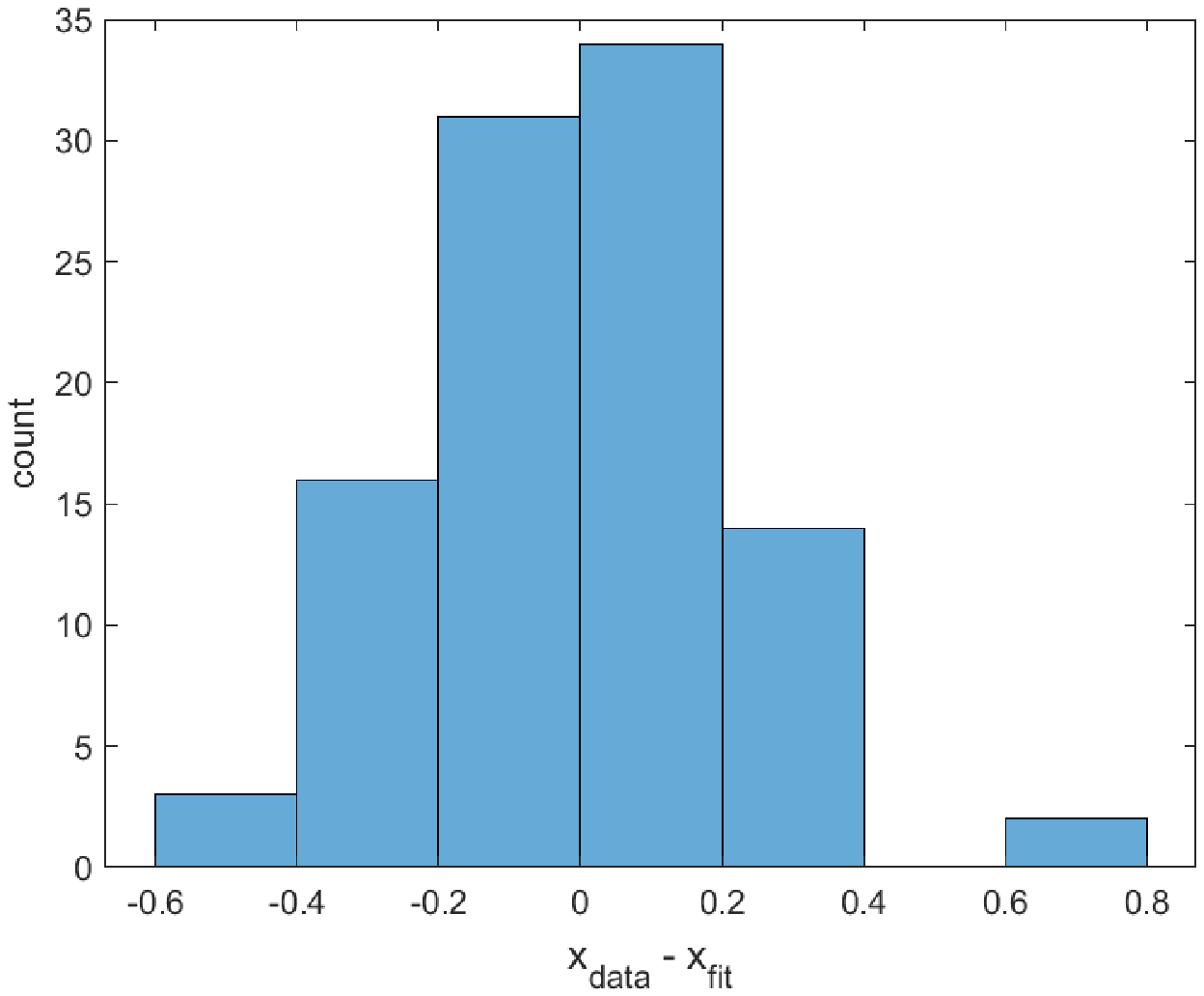}		
		\label{fig:para_est1_hist_lmx_alpha0p3beta0p5}
	\end{subfigure}
	\hfill
	\begin{subfigure}[b]{0.48\textwidth}
		\centering
		\includegraphics[width=1\textwidth]{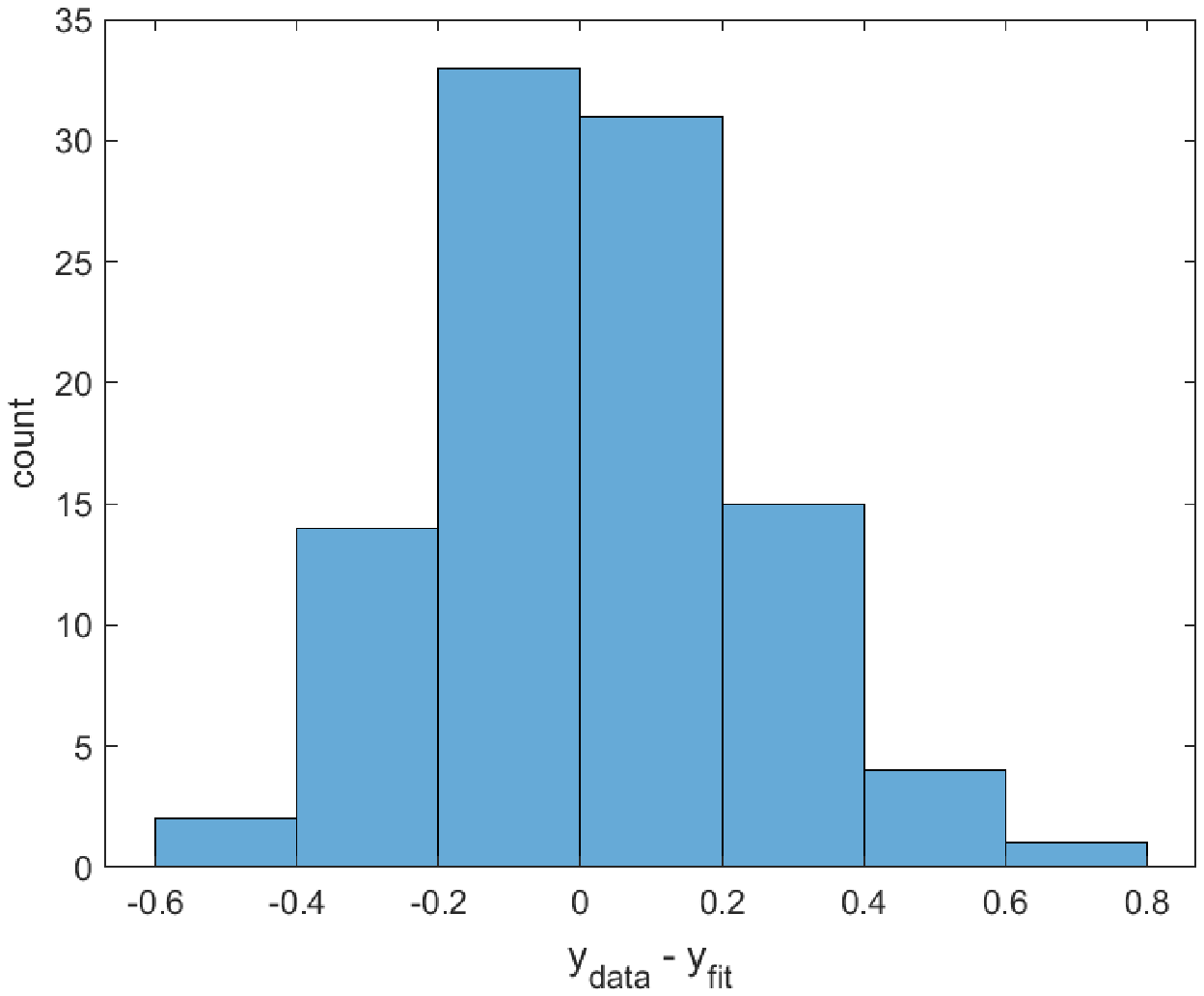}		
		\label{fig:para_est1_hist_lmy_alpha0p3beta0p5}
	\end{subfigure}
	\caption{Histogram of the errors between the data and curve fit with Levenberg–Marquardt algorithm for Example 1}
	\label{fig:para_est1_hist_lm_alpha0p3beta0p5}
\end{figure}

The fitted parameters are off by about 0.42\% in $\alpha$ and 0.05\% in $\beta.$ We have good recovery of the original parameters. Figure \ref{fig:para_est1_hist_lm_alpha0p3beta0p5} shows the histogram of the difference between the data values and the best-fit.

\subsection*{Computaion of Estimates with Trust-region algorithm}

\begin{figure}[H]
	\centering
	\begin{subfigure}[b]{0.48\textwidth}
		\centering
		\includegraphics[width=1\textwidth]{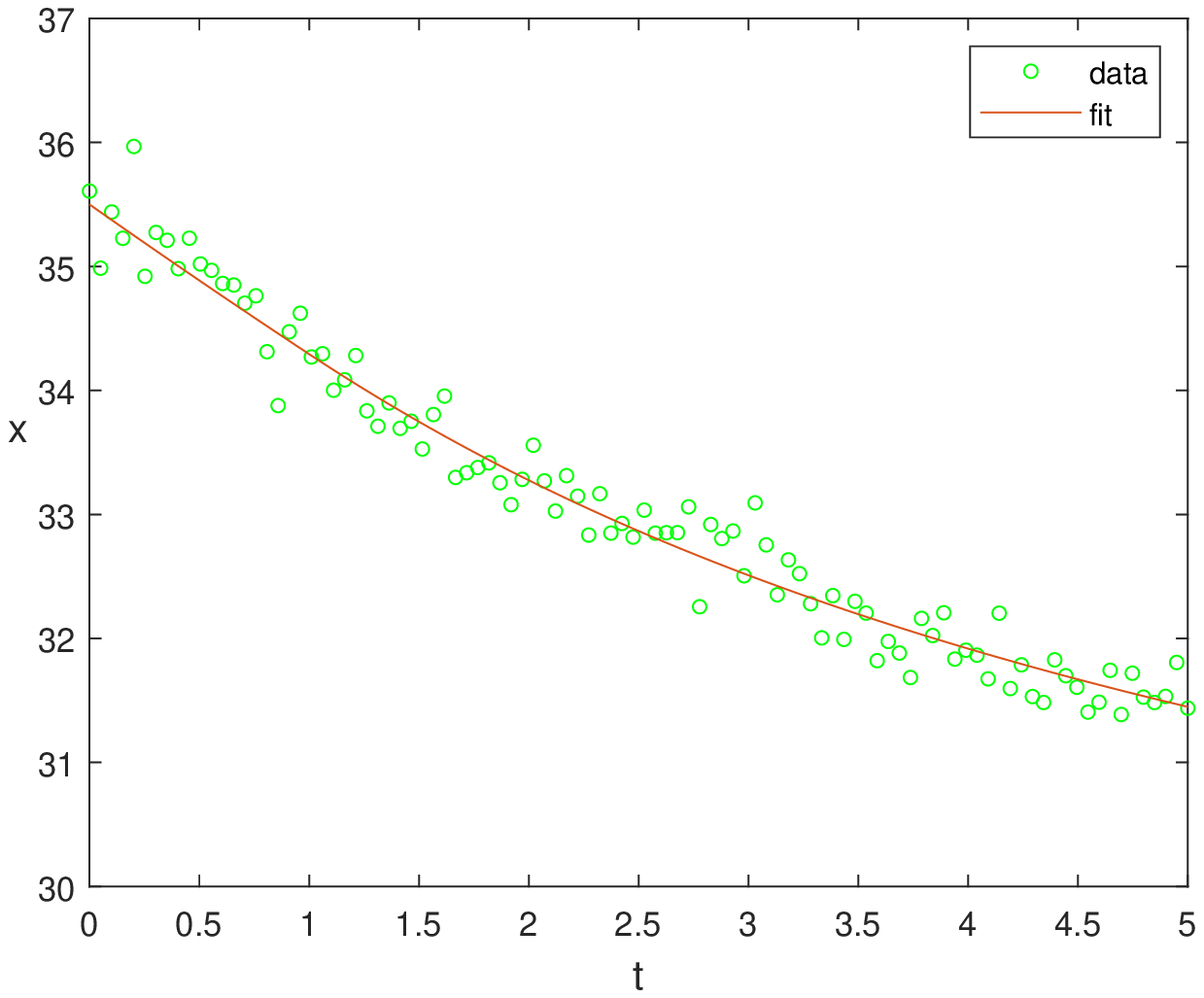}		
		\label{fig:para_est1_trx_alpha0p3beta0p5}
	\end{subfigure}
	\hfill
	\begin{subfigure}[b]{0.48\textwidth}
		\centering
		\includegraphics[width=1\textwidth]{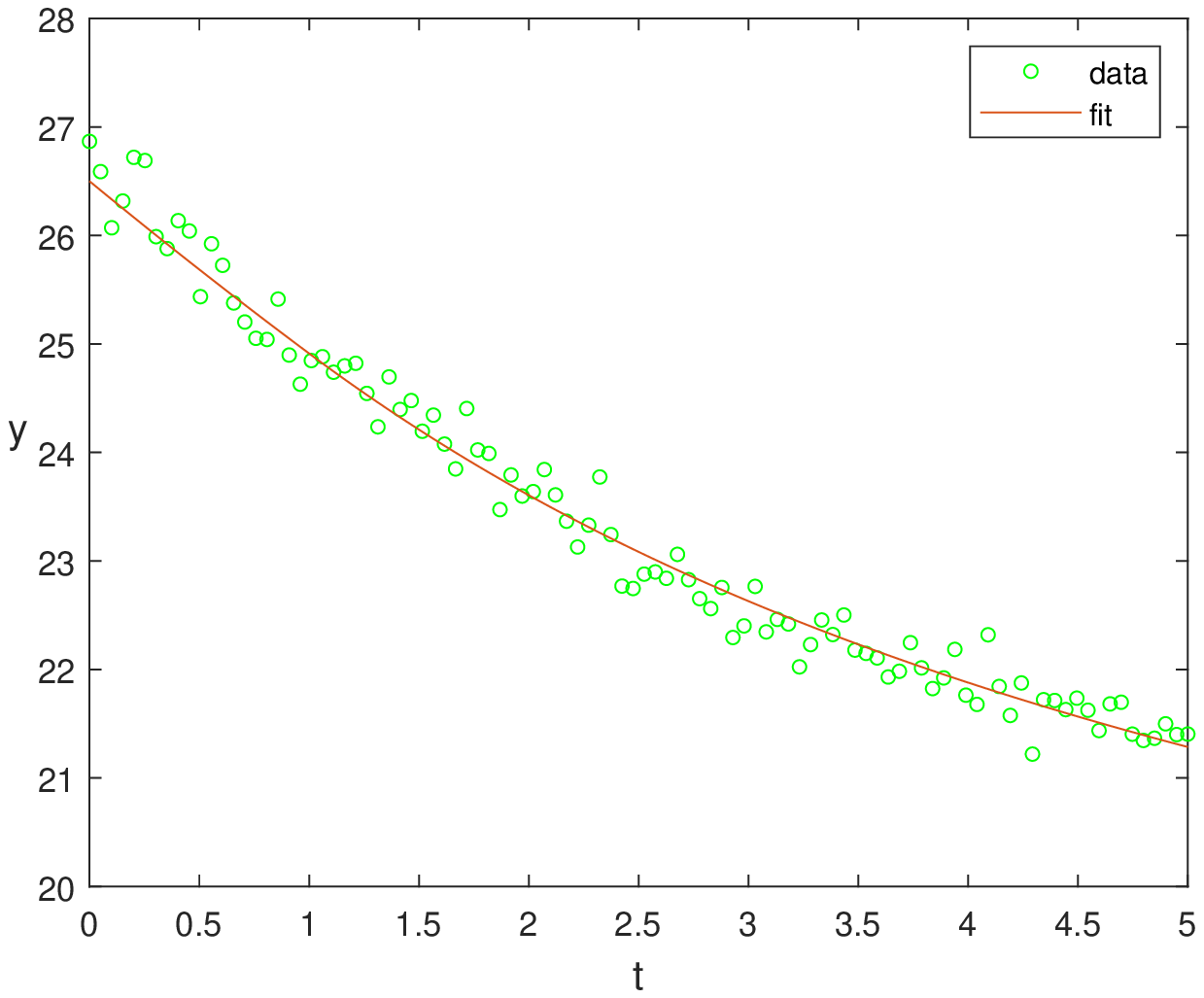}		
		\label{fig:para_est1_try_alpha0p3beta0p5}
	\end{subfigure}
	\caption{Data and curve fit for $x$ and $y$ with Trust-region algorithm for Example 1}
	\label{fig:para_est1_tr_alpha0p3beta0p5}
\end{figure}

\begin{table}[H]
	\begin{center}
		\caption{Estimation of $\alpha$ and $\beta$ with  Trust-region algorithm for Example 1}
		\label{table:ex1_tr}
		\begin{tabular}{ccccr}
			\hline
			Iteration &  \multicolumn{1}{p{2cm}}{\centering Function\\ count} & Residual & Norm of step &  \multicolumn{1}{p{2cm}}{\centering First-order\\ optimality}      \\
			\hline
			0 & 3  & 738.83  &             & 1.96e+04    \\
			1 & 6  & 17.6464 & 0.102945    & 1.59e+03 \\
			2 & 9  & 9.39188 & 0.0139484   & 17.9     \\
			3 & 12 & 9.39033 & 0.00020538  & 0.00338  \\
			4 & 15 & 9.39033 & 6.93959e-08 & 1.28e-05	\\	
			\hline
		\end{tabular}
	\end{center}
\end{table}

The data points, the curve fit are plotted in Figure \ref{fig:para_est1_tr_alpha0p3beta0p5}.  Table \ref{table:ex1_tr} contains the values of  iteration count, function count, residual, norm of the step and first-order optimality. We used the default settings of   the Trust-region algorithm in Matlab. The optimization stopped because the relative sum of squares  is changing by less than function tolerance of 1.000000e-06. After four iterations it gives the best fit of the parameters.

\begin{table}[H]
	\begin{center}
		\caption{True and best fit of  $\alpha$ and $\beta$ for Example 1}
		\label{table:ex1_true_fit_tr}
		\begin{tabular}{|c|c|c|c|}
			\hline
			&Initial & True & Best fit \\
			\hline
			$\alpha$ & 0.3 & 0.5  & 0.5021  \\ \hline
			$\beta$  & 0.5 & 0.8 & 0.7996 \\
			\hline
		\end{tabular}
	\end{center}
\end{table}

\begin{figure}[H]
	\centering
	\begin{subfigure}[b]{0.48\textwidth}
		\centering
		\includegraphics[width=1\textwidth]{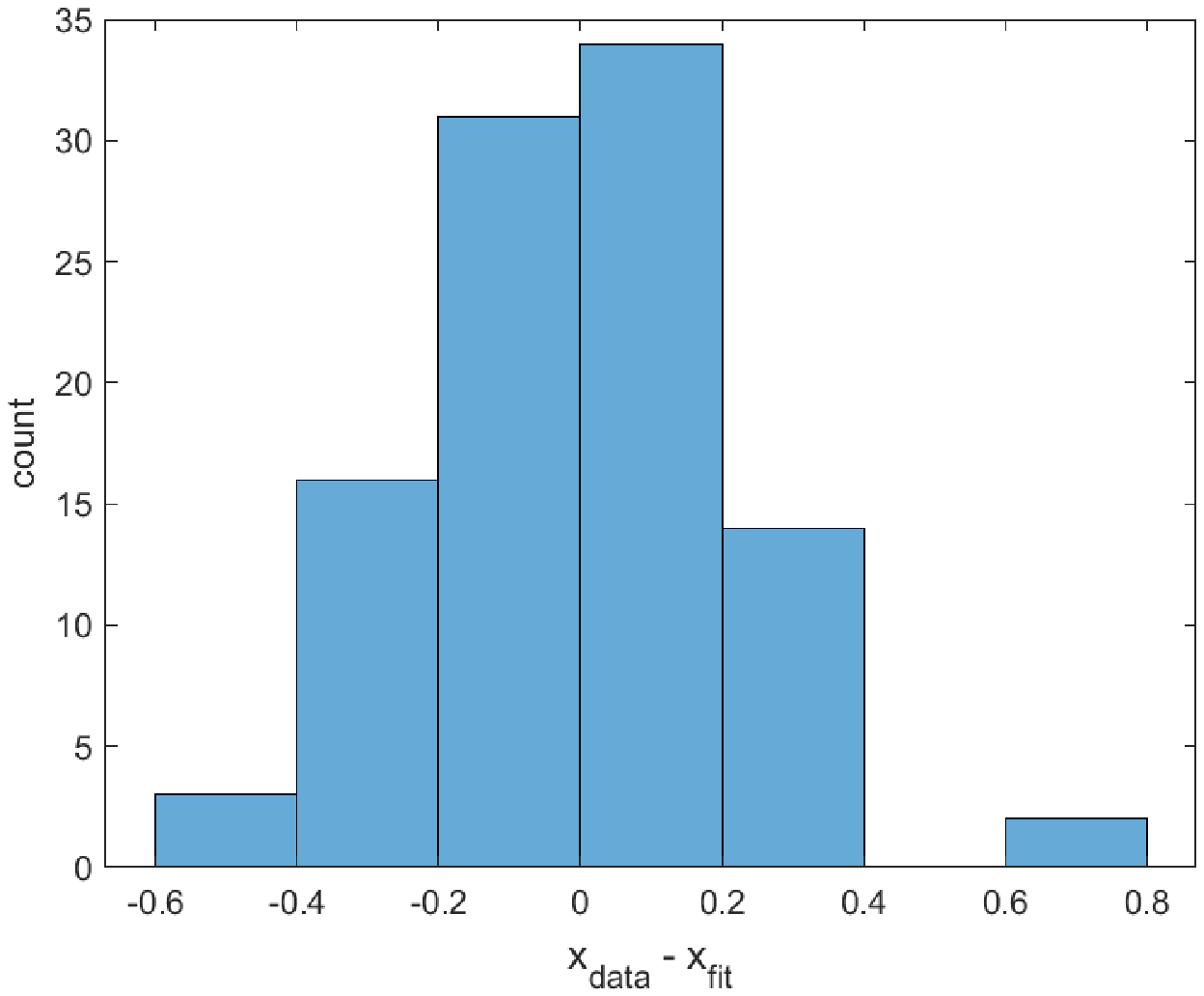}		
		\label{fig:para_est1_hist_trx_alpha0p3beta0p5}
	\end{subfigure}
	\hfill
	\begin{subfigure}[b]{0.48\textwidth}
		\centering
		\includegraphics[width=1\textwidth]{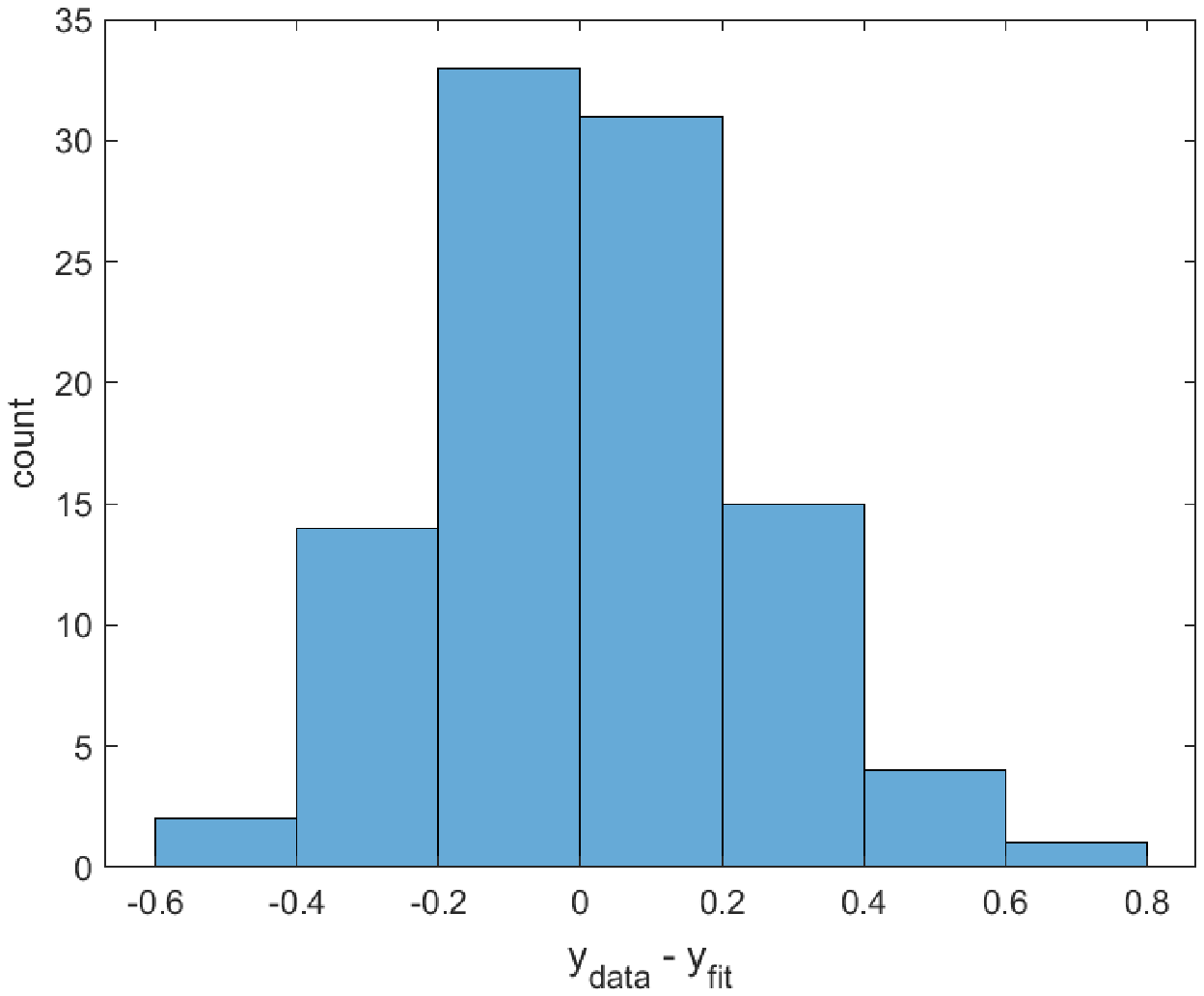}		
		\label{fig:para_est1_hist_try_alpha0p3beta0p5}
	\end{subfigure}
	\caption{Histogram of the errors between the data and curve fit with Trust-region algorithm for Example 1}
	\label{fig:para_est1_hist_tr_alpha0p3beta0p5}
\end{figure}

The fitted parameters are off by about 0.42\% in $\alpha$ and 0.05\% in $\beta$ which is the same as in Levenberg–Marquardt algorithm. We have good recovery again of the original parameters. Figure \ref{fig:para_est1_hist_tr_alpha0p3beta0p5} shows the histogram of the difference between the data values and the best-fit. 

The two algorithms found the same solution taking the same number of iterations. 

\section*{Example 2}

In this example we just change the initial starting point of the parameters. We use $\alpha = 0.01$ and $\beta = 0.01.$
\subsection*{Computation of Estimates with Levenberg–Marquardt algorithm}

\begin{figure}[H]
	\centering
	\begin{subfigure}[b]{0.48\textwidth}
		\centering
		\includegraphics[width=1\textwidth]{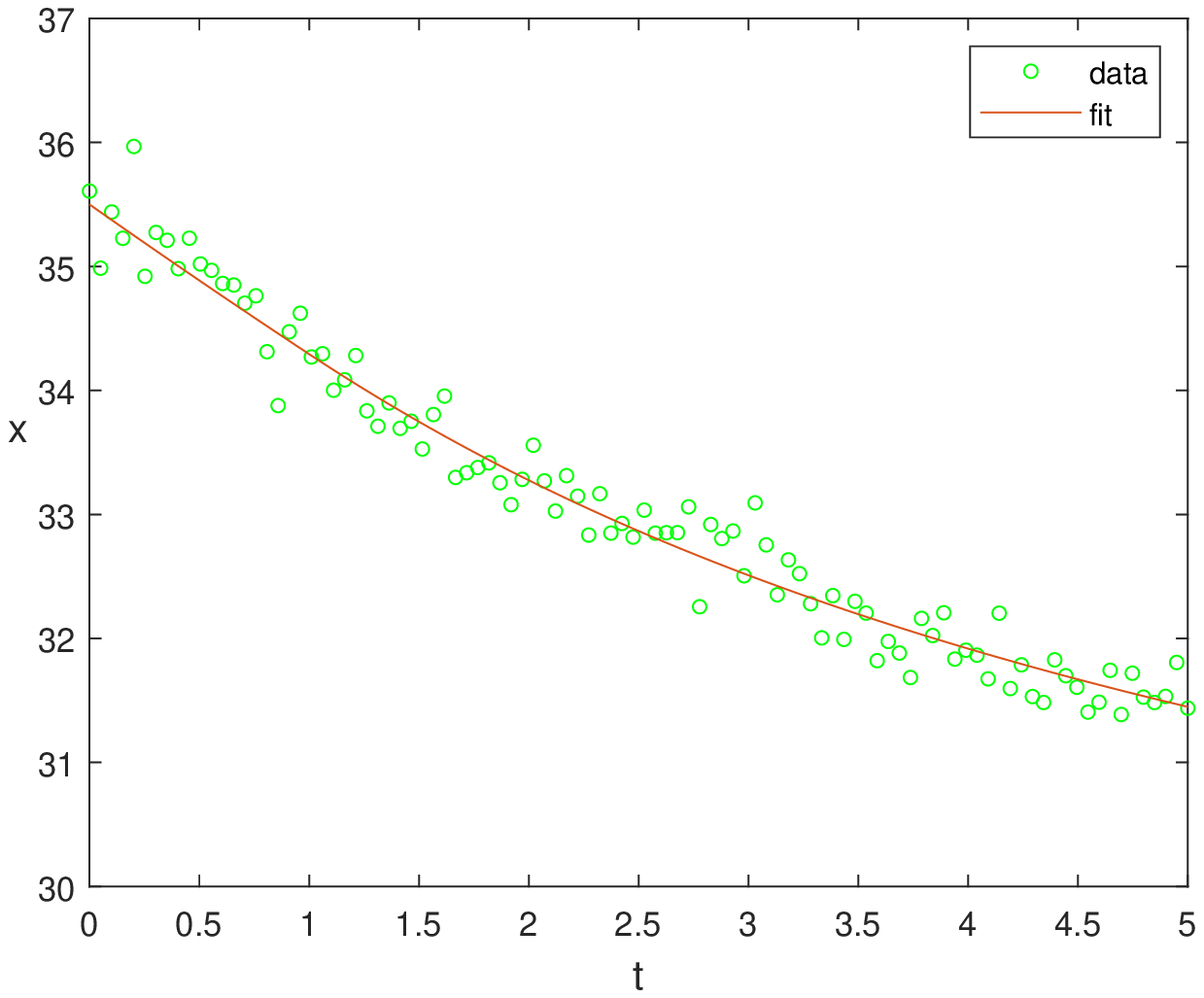}		
		\label{fig:para_est1_lmx_alpha0p01beta0p01}
	\end{subfigure}
	\hfill
	\begin{subfigure}[b]{0.48\textwidth}
		\centering
		\includegraphics[width=1\textwidth]{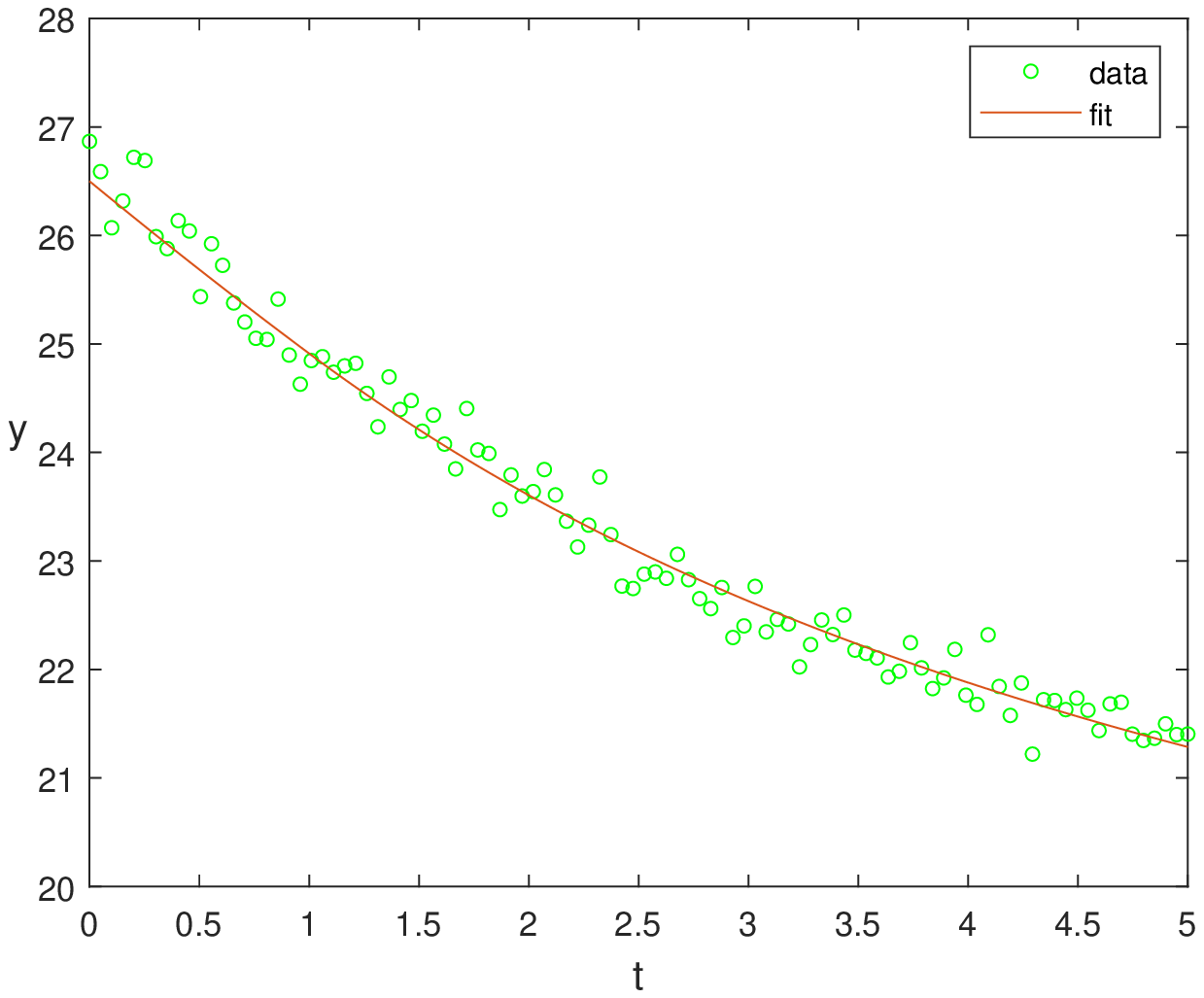}		
		\label{fig:para_est1_lmy_alpha0p01beta0p01}
	\end{subfigure}
	\caption{Data and curve fit for $x$ and $y$ with Levenberg–Marquardt algorithm for Example 2}
	\label{fig:para_est1_lm_alpha0p01beta0p01}
\end{figure}

\begin{table}[H]
	\begin{center}
		\caption{Estimation of $\alpha$ and $\beta$ with  Levenberg–Marquardt algorithm for Example 2}
		\label{table:ex2_lm}
		\begin{tabular}{cccccr}
			\hline
			Iteration &  \multicolumn{1}{p{2cm}}{\centering Function\\ count} & Residual &  \multicolumn{1}{p{2cm}}{\centering First-order\\ optimality} & Lambda & Norm of step    \\
			\hline
			0 & 3  & 6862.57 & 7.94e+03 & 0.01   &           \\
			1 & 6  & 346.729 & 1.08e+03 & 0.001  & 0.667034    \\
			2 & 9  & 10.8027 & 60.7     & 0.0001 & 0.245493    \\
			3 & 12 & 9.39036 & 0.32     & 1e-05  & 0.0187851   \\
			4 & 15 & 9.39033 & 7.49e-05 & 1e-06  & 9.29755e-05 \\	
			5 & 18 & 9.39033 & 1.62e-06 & 1e-07  & 8.75649e-09 \\	
			\hline
		\end{tabular}
	\end{center}
\end{table}

 After five iterations it gives the same best fit of the parameters as in Example 1.

\begin{table}[H]
	\begin{center}
		\caption{True and best fit of  $\alpha$ and $\beta$ for Example 2}
		\label{table:ex2_true_fit_lm}
		\begin{tabular}{|c|c|c|c|}
			\hline
			&Initial & True & Best fit \\
			\hline
			$\alpha$ & 0.01 & 0.5  & 0.5021  \\ \hline
			$\beta$  & 0.01 & 0.8 & 0.7996 \\
			\hline
		\end{tabular}
	\end{center}
\end{table}

\subsection*{Computation of Estimates with Trust-region algorithm}

\begin{figure}[H]
	\centering
	\begin{subfigure}[b]{0.48\textwidth}
		\centering
		\includegraphics[width=1\textwidth]{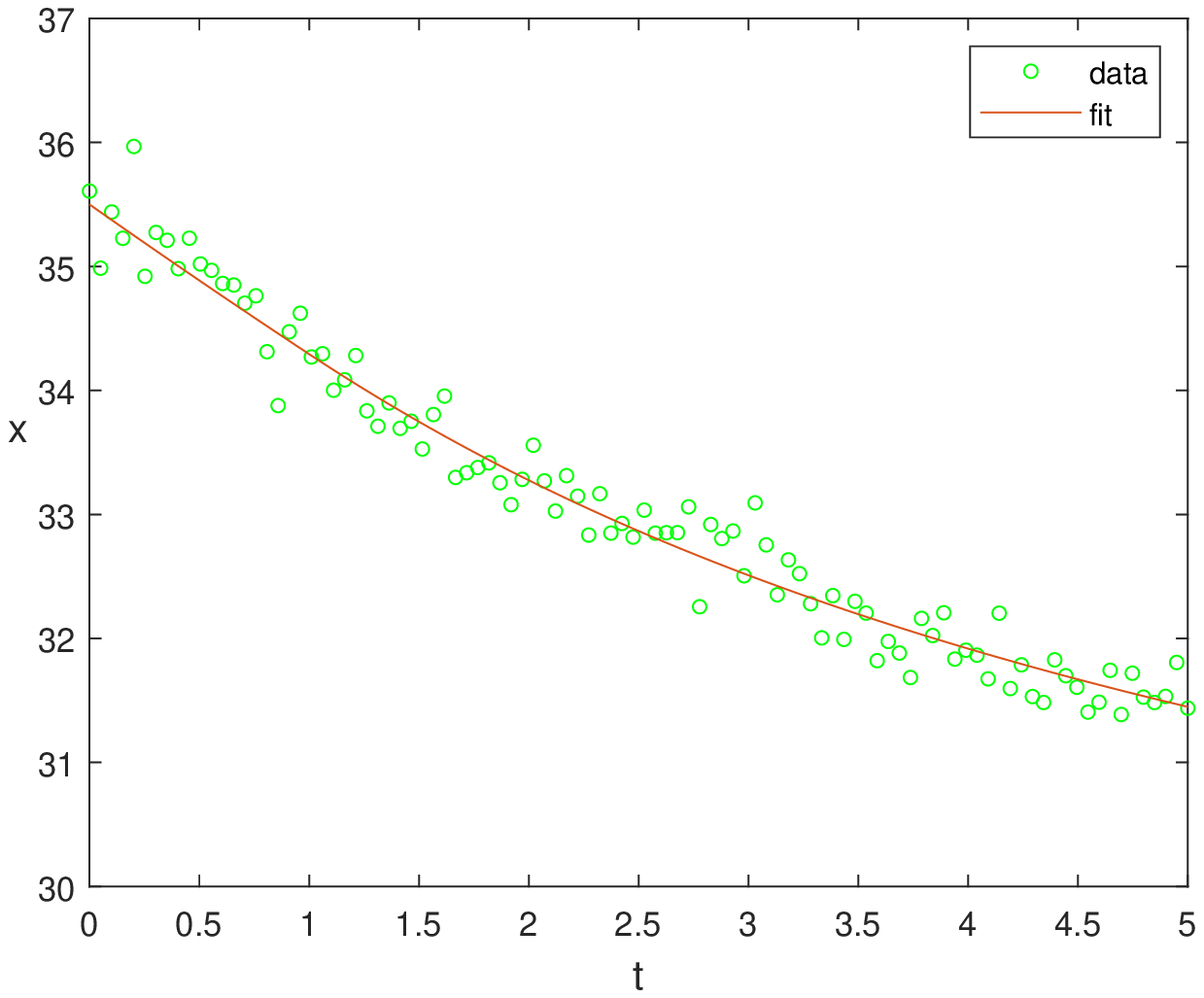}		
		\label{fig:para_est1_trx_alpha0p01beta0p01}
	\end{subfigure}
	\hfill
	\begin{subfigure}[b]{0.48\textwidth}
		\centering
		\includegraphics[width=1\textwidth]{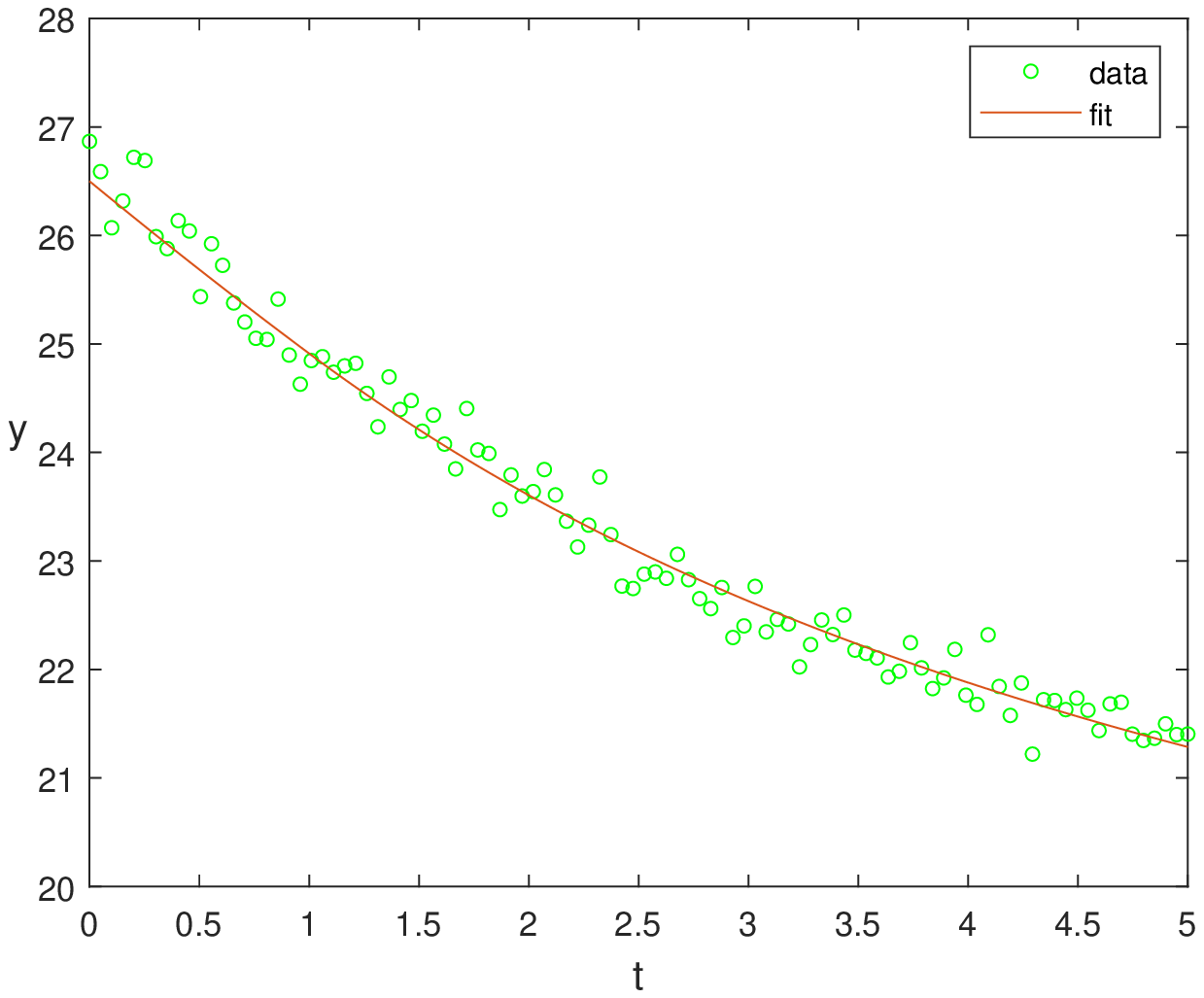}		
		\label{fig:para_est1_try_alpha0p01beta0p01}
	\end{subfigure}
	\caption{Data and curve fit for $x$ and $y$ with Trust-region algorithm for Example 2}
	\label{fig:para_est1_tr_alpha0p01beta0p01}
\end{figure}

\begin{table}[H]
	\begin{center}
		\caption{Estimation of $\alpha$ and $\beta$ with  Trust-region algorithm for Example 2}
		\label{table:ex2_tr}
		\begin{tabular}{ccccr}
			\hline
			Iteration &  \multicolumn{1}{p{2cm}}{\centering Function\\ count} & Residual & Norm of step &  \multicolumn{1}{p{2cm}}{\centering First-order\\ optimality}      \\
			\hline
			0 & 3  & 41238.5 & 1.61e+04    &          \\
			1 & 6  & 9601.81 & 1.95587     & 6.79e+03 \\
			2 & 9  & 1550.04 & 0.941773    & 2.05e+03 \\
			3 & 12 & 144.481 & 0.451259    & 483      \\
			4 & 15 & 12.7639 & 0.167543    & 69.4     \\
			5 & 18 & 9.39496 & 0.0306896   & 2.56     \\
			6 & 21 & 9.39033 & 0.00116804  & 0.00343  \\
			7 & 24 & 9.39033 & 1.51337e-06 & 8.85e-06 \\
			\hline
		\end{tabular}
	\end{center}
\end{table}

 After seven iterations it gives the best fit of the parameters.

\begin{table}[H]
	\begin{center}
		\caption{True and best fit of  $\alpha$ and $\beta$ for Example 2}
		\label{table:ex2_true_fit_tr}
		\begin{tabular}{|c|c|c|c|}
			\hline
			&Initial & True & Best fit \\
			\hline
			$\alpha$ & 0.3 & 0.5  & 0.5021  \\ \hline
			$\beta$  & 0.5 & 0.8 & 0.7996 \\
			\hline
		\end{tabular}
	\end{center}
\end{table}

\section*{Example 3}

In this example we  repeat example 1 but the  noise in our data is increased. The noise  has a normal distribution with a mean of zero and a standard deviation of 0.40. 
\begin{equation}
	\begin{aligned}
		x_i &= x(t_i; \alpha = 0.5, \beta = 0.8) + \mathcal{N}(0, 0.40)\\	
		y_i &= y(t_i; \alpha = 0.5, \beta = 0.8) + \mathcal{N}(0, 0.40)	
	\end{aligned}
\end{equation}

\subsection*{Computation of Estimates with Levenberg–Marquardt algorithm}

\begin{figure}[H]
	\centering
	\begin{subfigure}[b]{0.48\textwidth}
		\centering
		\includegraphics[width=1\textwidth]{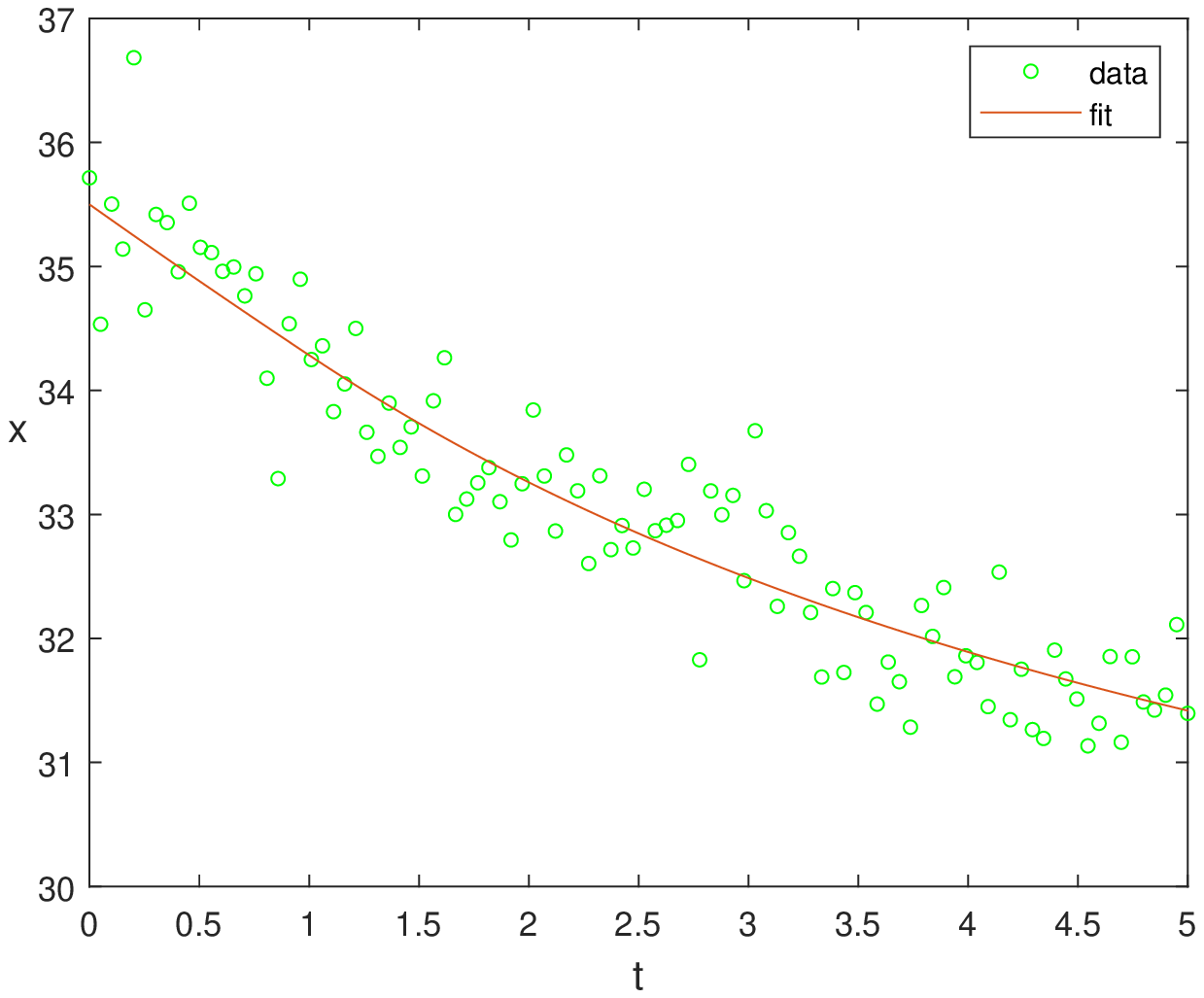}		
		\label{fig:para_est3_lmx_alpha0p3beta0p5}
	\end{subfigure}
	\hfill
	\begin{subfigure}[b]{0.48\textwidth}
		\centering
		\includegraphics[width=1\textwidth]{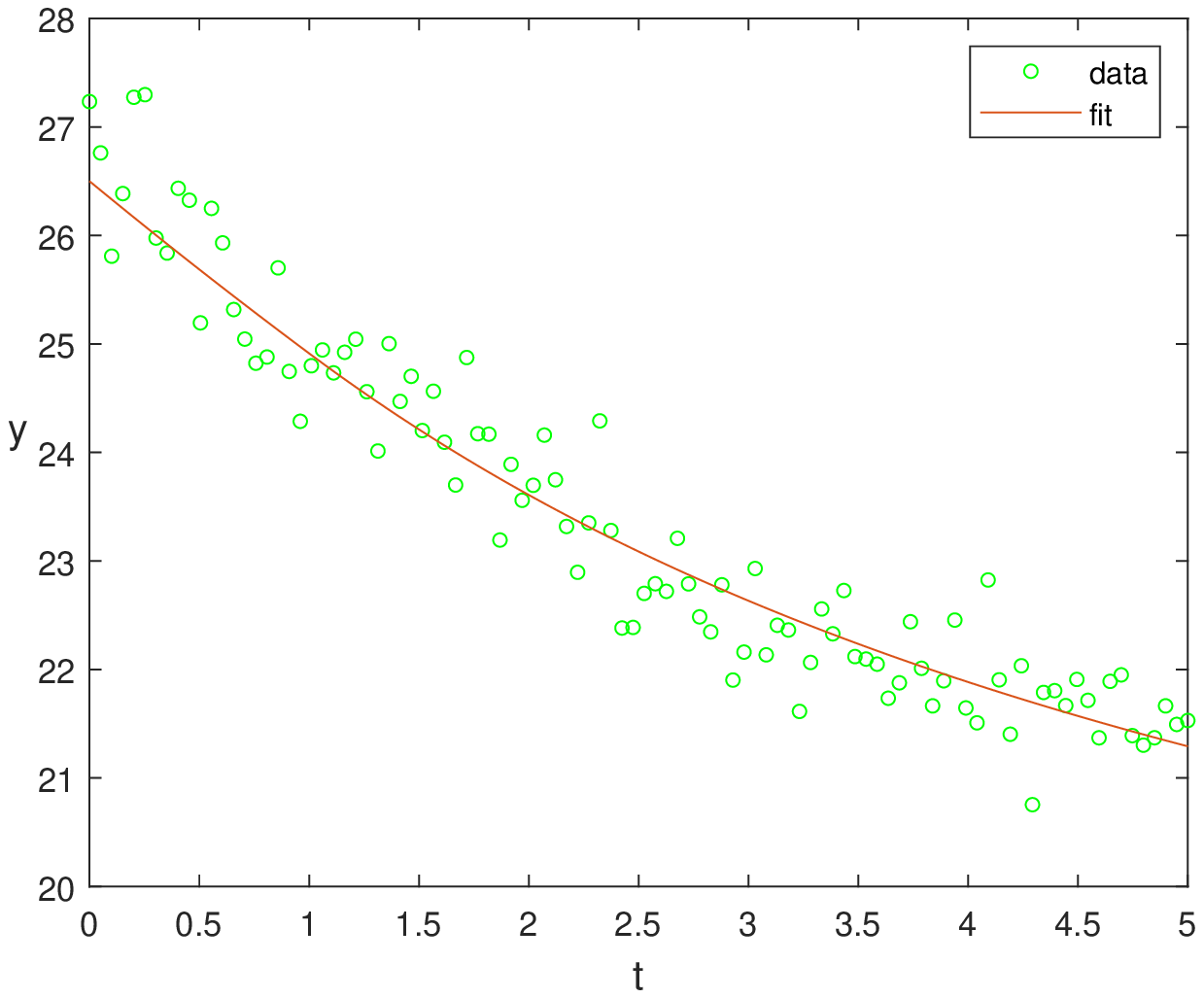}		
		\label{fig:para_est3_lmy_alpha0p3beta0p5}
	\end{subfigure}
	\caption{Data and curve fit for $x$ and $y$ with Levenberg–Marquardt algorithm for Example 3}
	\label{fig:para_est3_lm_alpha0p3beta0p5}
\end{figure}

\begin{table}[H]
	\begin{center}
		\caption{Estimation of $\alpha$ and $\beta$ with  Levenberg–Marquardt algorithm for Example 3}
		\label{table:ex3_lm}
		\begin{tabular}{cccccr}
			\hline
			Iteration &  \multicolumn{1}{p{2cm}}{\centering Function\\ count} & Residual &  \multicolumn{1}{p{2cm}}{\centering First-order\\ optimality} & Lambda & Norm of step    \\
			\hline
			0 & 3  & 773.538 & 2.04e+03 & 0.01   &             \\
			1 & 6  & 43.1368 & 141      & 0.001  & 0.327075    \\
			2 & 9  & 37.5617 & 0.957    & 0.0001 & 0.0349491   \\
			3 & 12 & 37.5613 & 0.000643 & 1e-05  & 0.000308975 \\
			4 & 15 & 37.5613 & 2.09e-06 & 1e-06  & 7.65485e-08 \\
			\hline
		\end{tabular}
	\end{center}
\end{table}

\begin{table}[H]
	\begin{center}
		\caption{True and best fit of  $\alpha$ and $\beta$ for Example 3}
		\label{table:ex3_true_fit_lm}
		\begin{tabular}{|c|c|c|c|}
			\hline
			&Initial & True & Best fit \\
			\hline
			$\alpha$ & 0.3 & 0.5  & 0.5042  \\ \hline
			$\beta$  & 0.5 & 0.8 & 0.7992 \\
			\hline
		\end{tabular}
	\end{center}
\end{table}

\begin{figure}[H]
	\centering
	\begin{subfigure}[b]{0.48\textwidth}
		\centering
		\includegraphics[width=1\textwidth]{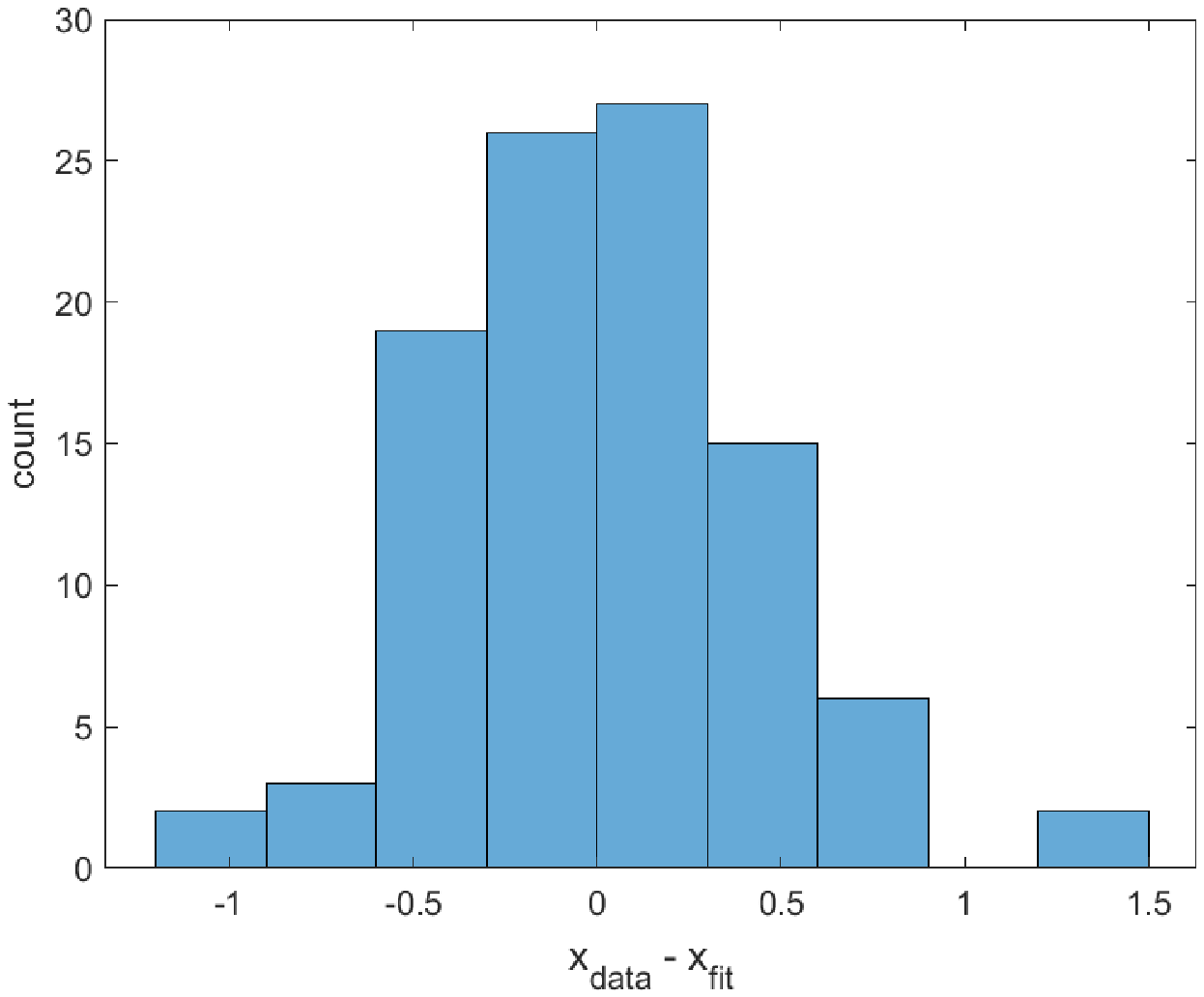}		
		\label{fig:para_est3_hist_lmx_alpha0p3beta0p5}
	\end{subfigure}
	\hfill
	\begin{subfigure}[b]{0.48\textwidth}
		\centering
		\includegraphics[width=1\textwidth]{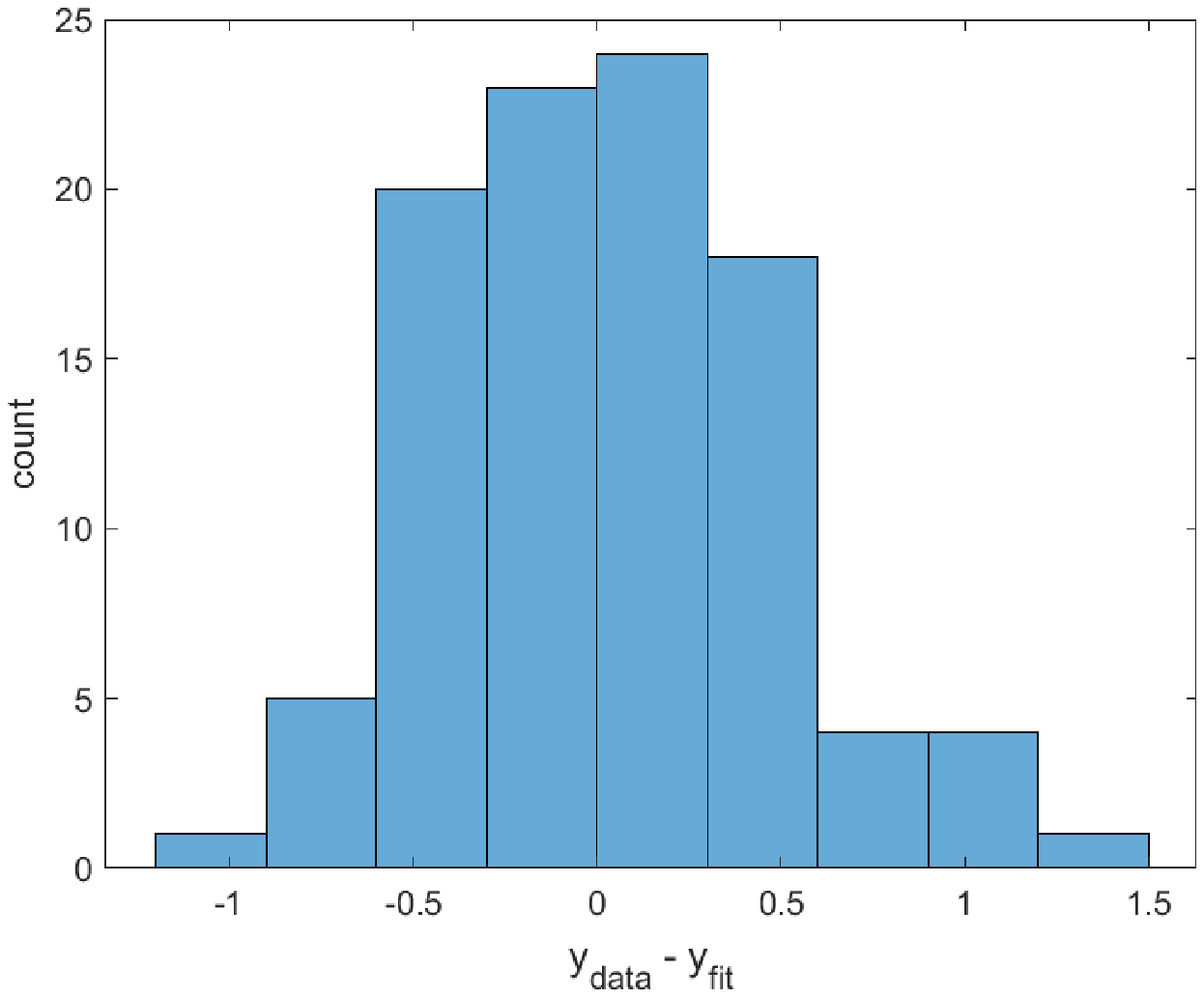}		
		\label{fig:para_est3_hist_lmy_alpha0p3beta0p5}
	\end{subfigure}
	\caption{Histogram of the errors between the data and curve fit with Levenberg–Marquardt algorithm for Example 3}
	\label{fig:para_est3_hist_lm_alpha0p3beta0p5}
\end{figure}

The fitted parameters are off by about 0.84\% in $\alpha$ and 0.1\% in $\beta.$  Figure \ref{fig:para_est3_hist_lm_alpha0p3beta0p5} shows the histogram of the difference between the data values and the best-fit.

\subsection*{Computaion of Estimates with Trust-region algorithm}

\begin{figure}[H]
	\centering
	\begin{subfigure}[b]{0.48\textwidth}
		\centering
		\includegraphics[width=1\textwidth]{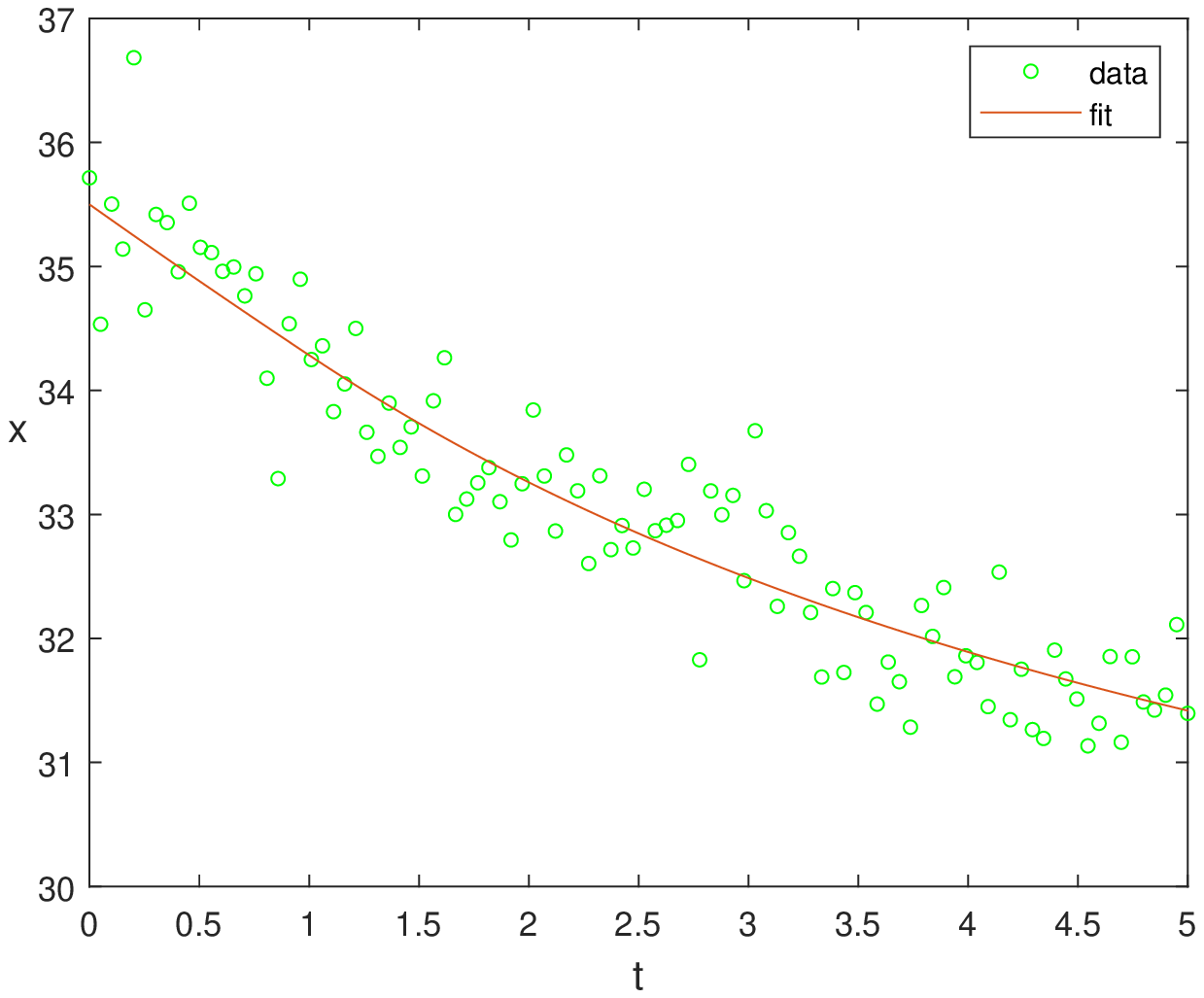}		
		\label{fig:para_est3_trx_alpha0p3beta0p5}
	\end{subfigure}
	\hfill
	\begin{subfigure}[b]{0.48\textwidth}
		\centering
		\includegraphics[width=1\textwidth]{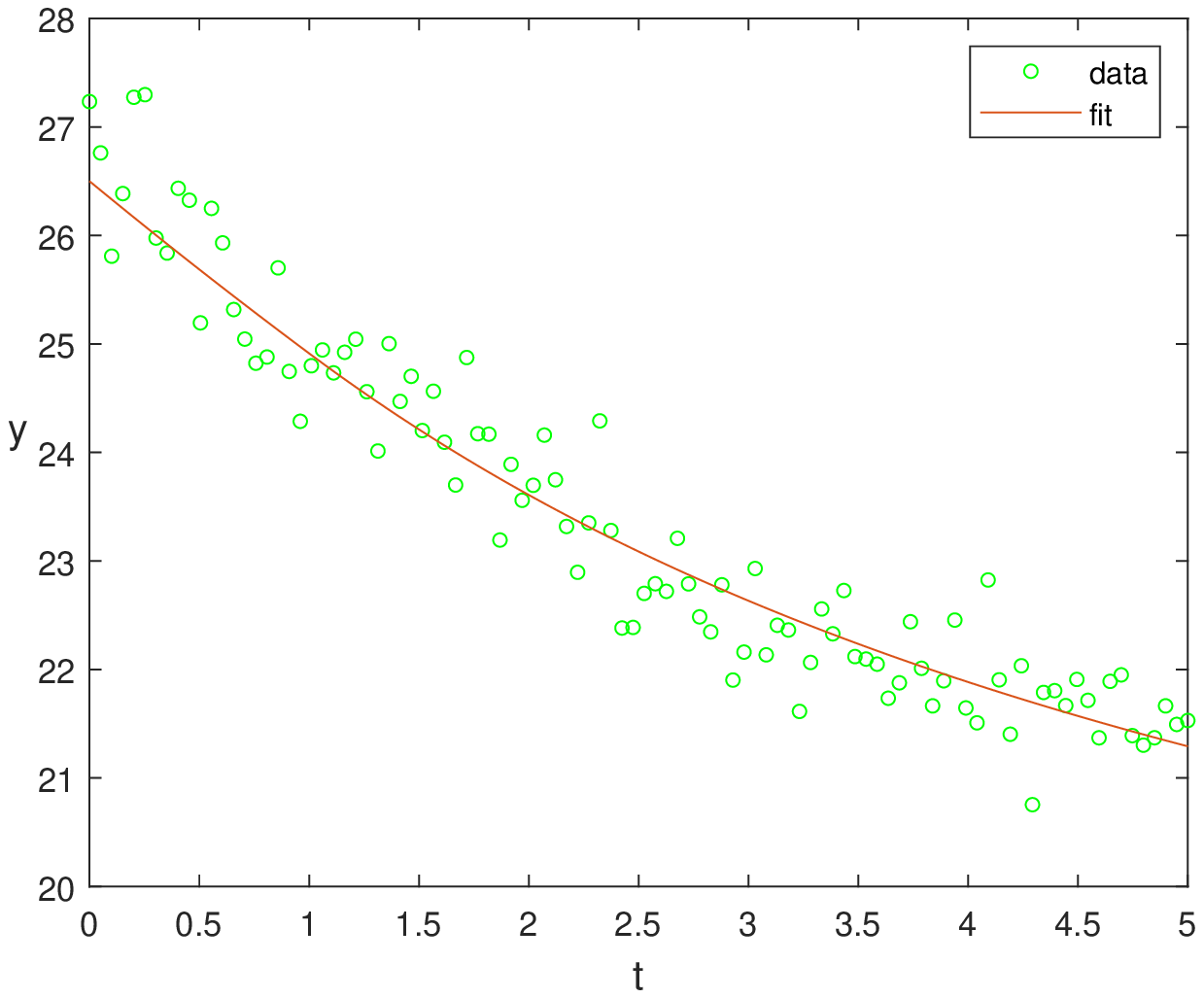}		
		\label{fig:para_est3_try_alpha0p3beta0p5}
	\end{subfigure}
	\caption{Data and curve fit for $x$ and $y$ with Trust-region algorithm for Example 3}
	\label{fig:para_est3_tr_alpha0p3beta0p5}
\end{figure}

\begin{table}[H]
	\begin{center}
		\caption{Estimation of $\alpha$ and $\beta$ with  Trust-region algorithm for Example 3}
		\label{table:ex3_tr}
		\begin{tabular}{ccccr}
			\hline
			Iteration &  \multicolumn{1}{p{2cm}}{\centering Function\\ count} & Residual & Norm of step &  \multicolumn{1}{p{2cm}}{\centering First-order\\ optimality}      \\
			\hline
			0 & 3  & 773.538 & 1.98e+04    &          \\
			1 & 6  & 45.902  & 0.10318     & 1.61e+03 \\
			2 & 9  & 37.5629 & 0.0139874   & 17.7     \\
			3 & 12 & 37.5613 & 0.000204213 & 0.00248  \\
			4 & 15 & 37.5613 & 1.87423e-07 & 3.07e-05 \\
			\hline
		\end{tabular}
	\end{center}
\end{table}

 After four iterations it gives the best fit of the parameters.

\begin{table}[H]
	\begin{center}
		\caption{True and best fit of  $\alpha$ and $\beta$ for Example 3}
		\label{table:ex3_true_fit_tr}
		\begin{tabular}{|c|c|c|c|}
			\hline
			&Initial & True & Best fit \\
			\hline
			$\alpha$ & 0.3 & 0.5  & 0.5042  \\ \hline
			$\beta$  & 0.5 & 0.8 & 0.7992 \\
			\hline
		\end{tabular}
	\end{center}
\end{table}

\begin{figure}[H]
	\centering
	\begin{subfigure}[b]{0.48\textwidth}
		\centering
		\includegraphics[width=1\textwidth]{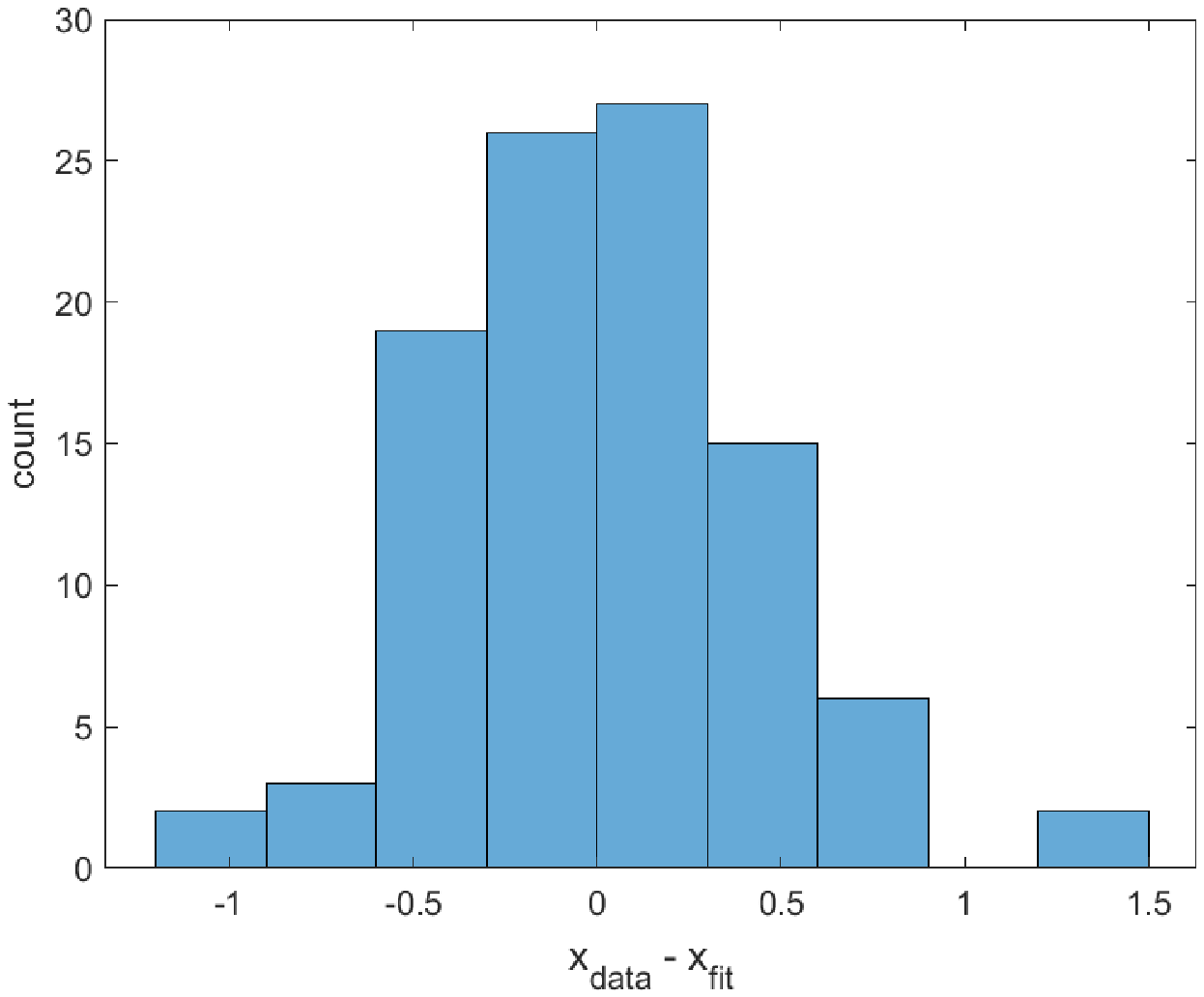}		
		\label{fig:para_est3_hist_trx_alpha0p3beta0p5}
	\end{subfigure}
	\hfill
	\begin{subfigure}[b]{0.48\textwidth}
		\centering
		\includegraphics[width=1\textwidth]{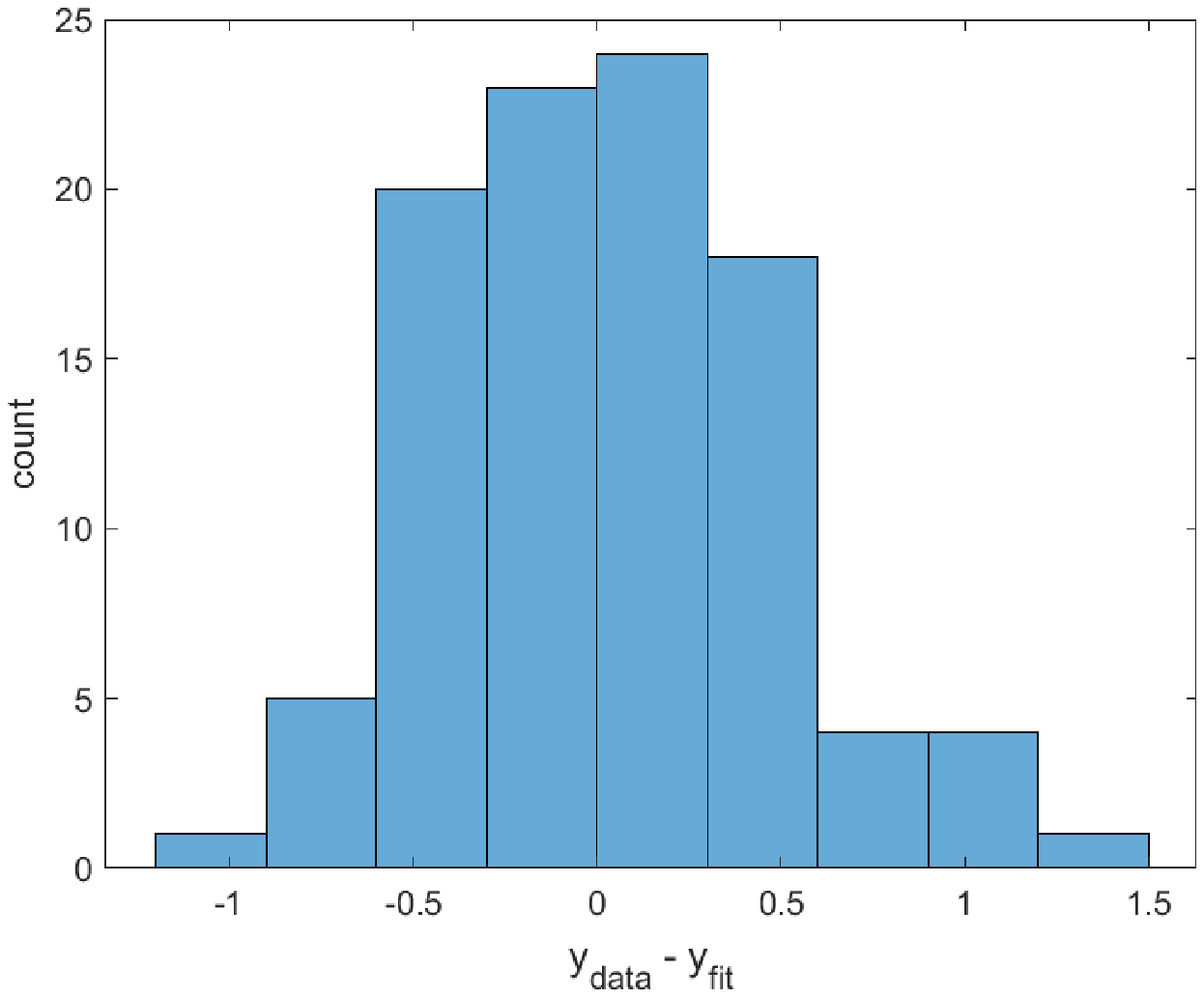}		
		\label{fig:para_est3_hist_try_alpha0p3beta0p5}
	\end{subfigure}
	\caption{Histogram of the errors between the data and curve fit with Trust-region algorithm for Example 3}
	\label{fig:para_est3_hist_tr_alpha0p3beta0p5}
\end{figure}

The fitted parameters are the same with both algorithms.   Figure \ref{fig:para_est3_hist_tr_alpha0p3beta0p5} shows the histogram of the difference between the data values and the best-fit. 

\section*{Example 4}

In this example we  repeat example 2 but the  noise in our data is increased. The noise  has a normal distribution with a mean of zero and a standard deviation of 0.40. 
\begin{equation}
	\begin{aligned}
		x_i &= x(t_i; \alpha = 0.5, \beta = 0.8) + \mathcal{N}(0, 0.40)\\	
		y_i &= y(t_i; \alpha = 0.5, \beta = 0.8) + \mathcal{N}(0, 0.40)	
	\end{aligned}
\end{equation}

\subsection*{Computation of Estimates with Levenberg–Marquardt algorithm}

\begin{figure}[H]
	\centering
	\begin{subfigure}[b]{0.48\textwidth}
		\centering
		\includegraphics[width=1\textwidth]{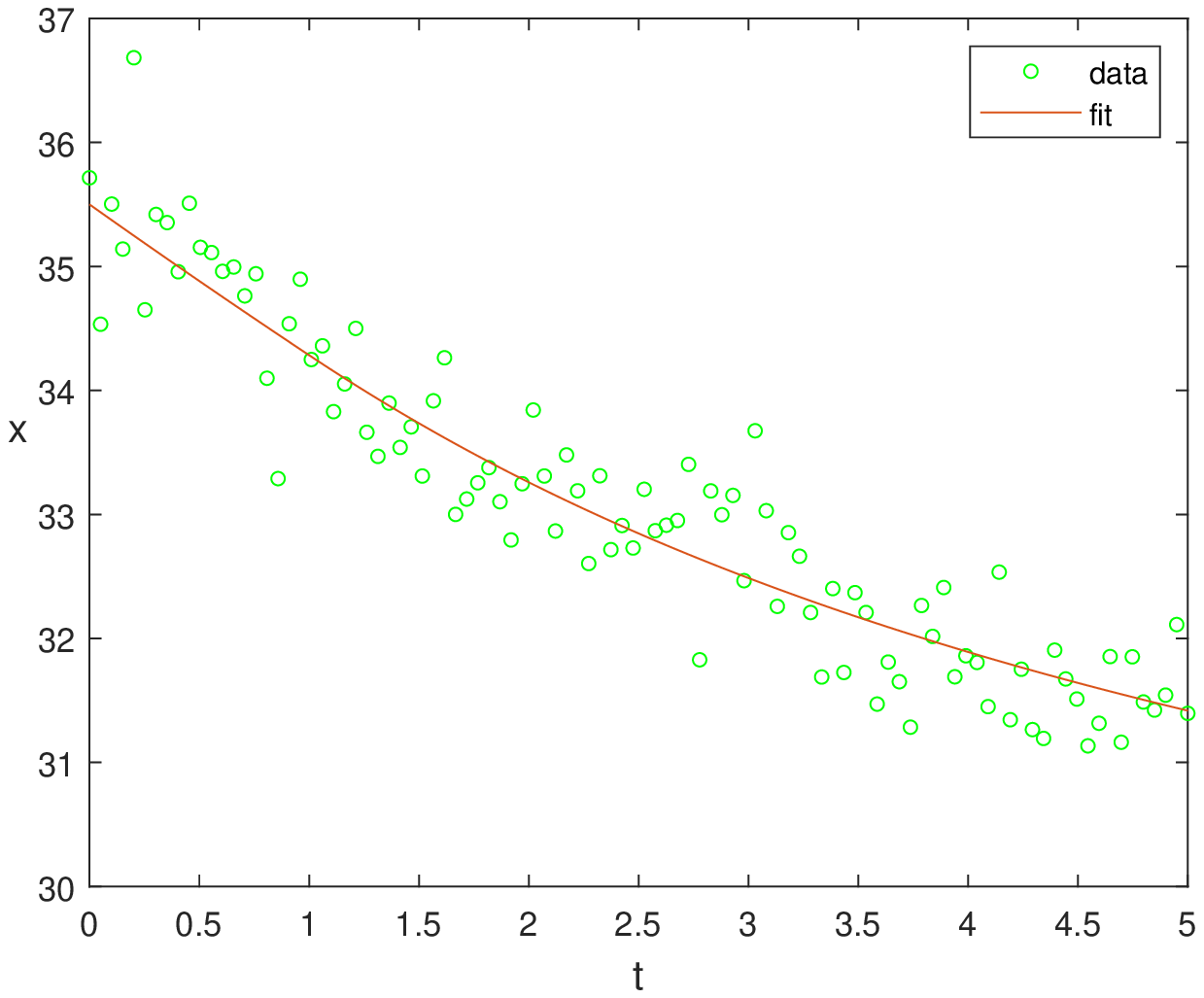}		
		\label{fig:para_est4_lmx_alpha0p01beta0p01}
	\end{subfigure}
	\hfill
	\begin{subfigure}[b]{0.48\textwidth}
		\centering
		\includegraphics[width=1\textwidth]{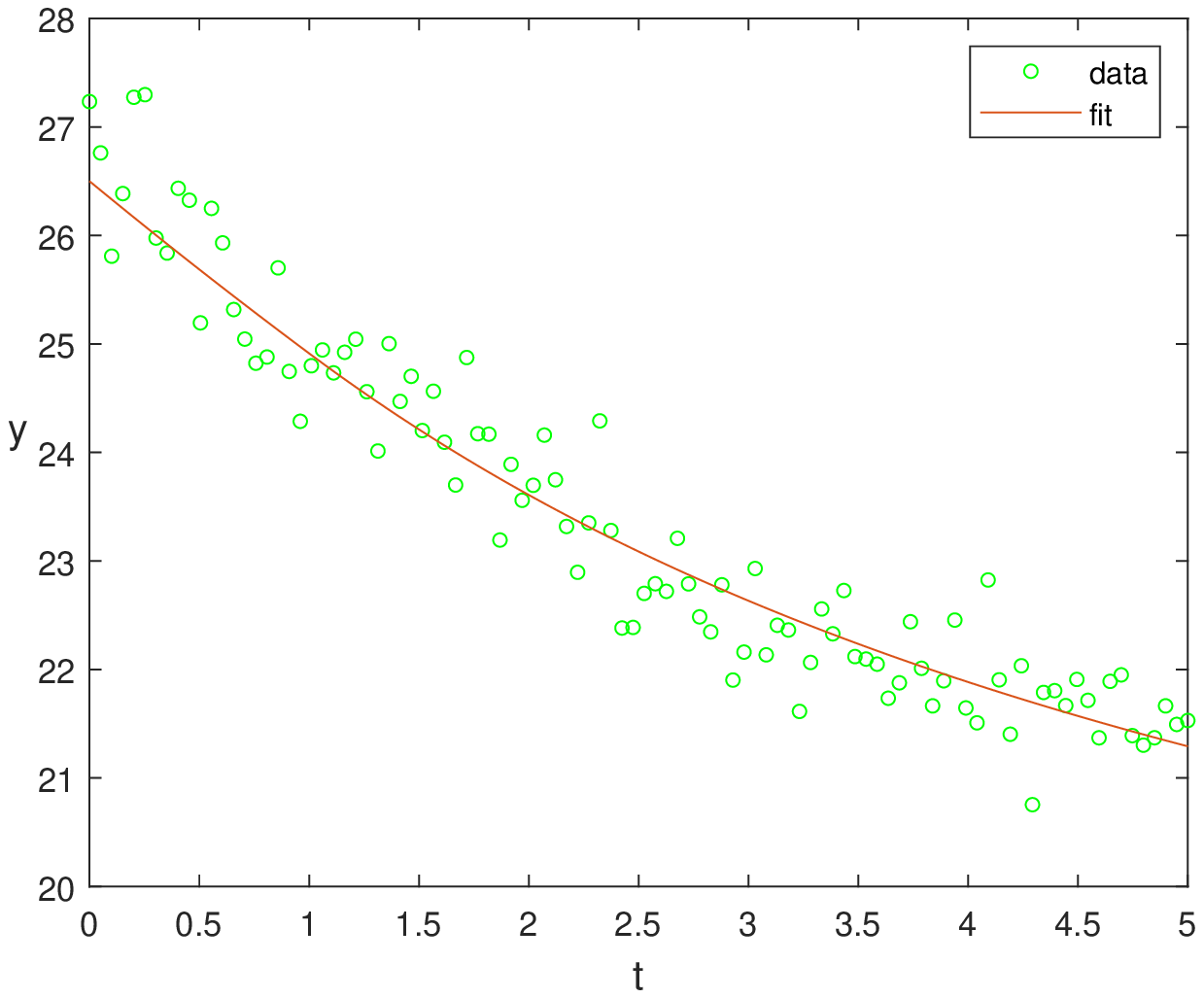}		
		\label{fig:para_est4_lmy_alpha0p01beta0p01}
	\end{subfigure}
	\caption{Data and curve fit for $x$ and $y$ with Levenberg–Marquardt algorithm for Example 4}
	\label{fig:para_est4_lm_alpha0p01beta0p01}
\end{figure}

\begin{table}[H]
	\begin{center}
		\caption{Estimation of $\alpha$ and $\beta$ with  Levenberg–Marquardt algorithm for Example 4}
		\label{table:ex4_lm}
		\begin{tabular}{cccccr}
			\hline
			Iteration &  \multicolumn{1}{p{2cm}}{\centering Function\\ count} & Residual &  \multicolumn{1}{p{2cm}}{\centering First-order\\ optimality} & Lambda & Norm of step    \\
			\hline
			0 & 3  & 6908.63 & 7.98e+03 & 0.01   &             \\
			1 & 6  & 376.107 & 1.08e+03 & 0.001  & 0.667467    \\
			2 & 9  & 38.9771 & 60.7     & 0.0001 & 0.245833    \\
			3 & 12 & 37.5614 & 0.318    & 1e-05  & 0.0188103   \\
			4 & 15 & 37.5613 & 0.000146 & 1e-06  & 9.14884e-05 \\
			\hline
		\end{tabular}
	\end{center}
\end{table}

 After four iterations it gives the same best fit of the parameters.

\begin{table}[H]
	\begin{center}
		\caption{True and best fit of  $\alpha$ and $\beta$ for Example 4}
		\label{table:ex4_true_fit_lm}
		\begin{tabular}{|c|c|c|c|}
			\hline
			&Initial & True & Best fit \\
			\hline
			$\alpha$ & 0.01 & 0.5  & 0.5042  \\ \hline
			$\beta$  & 0.01 & 0.8 & 0.7992 \\
			\hline
		\end{tabular}
	\end{center}
\end{table}

\subsection*{Computation of Estimates with Trust-region algorithm}

\begin{figure}[H]
	\centering
	\begin{subfigure}[b]{0.48\textwidth}
		\centering
		\includegraphics[width=1\textwidth]{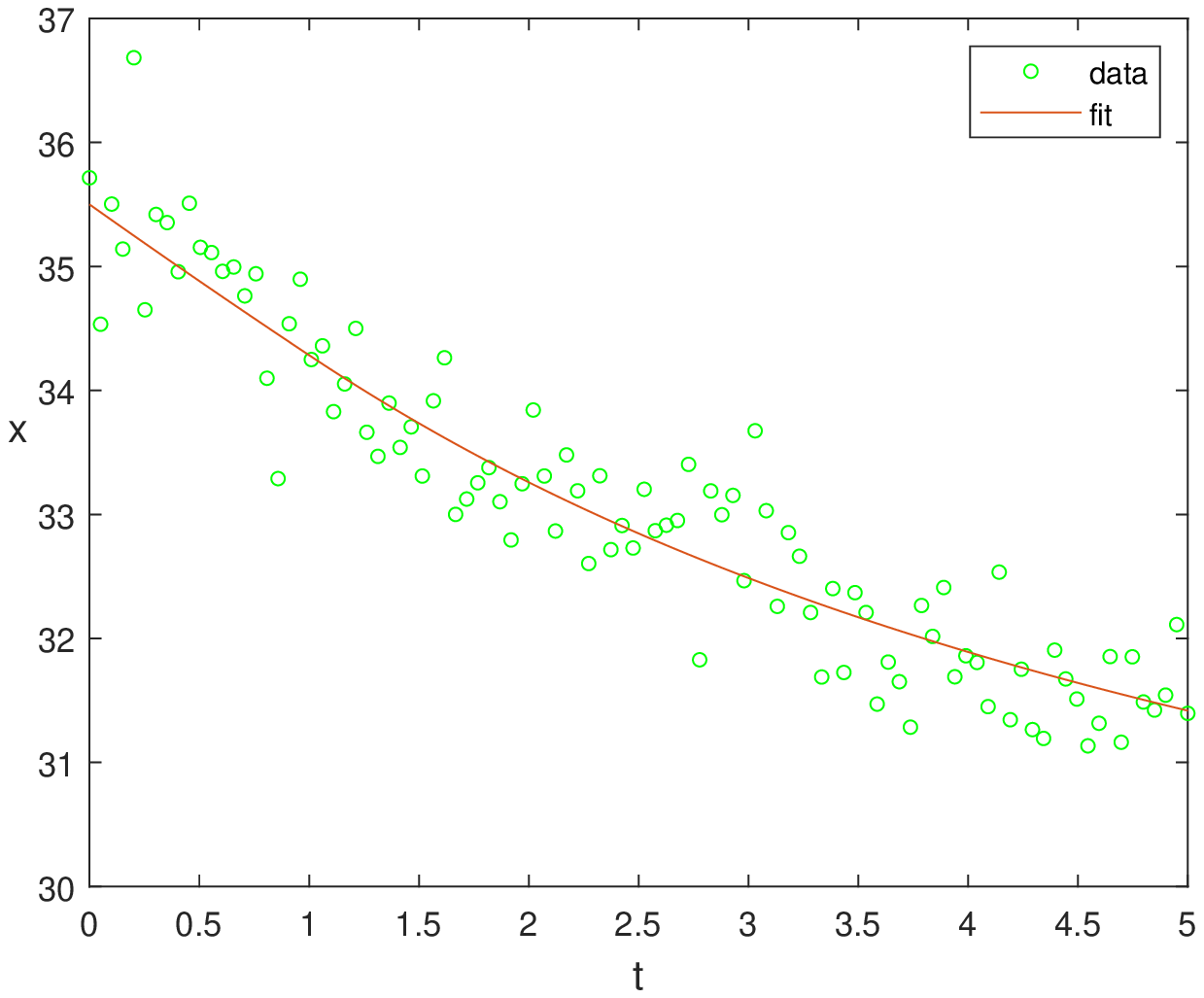}		
		\label{fig:para_est4_trx_alpha0p01beta0p01}
	\end{subfigure}
	\hfill
	\begin{subfigure}[b]{0.48\textwidth}
		\centering
		\includegraphics[width=1\textwidth]{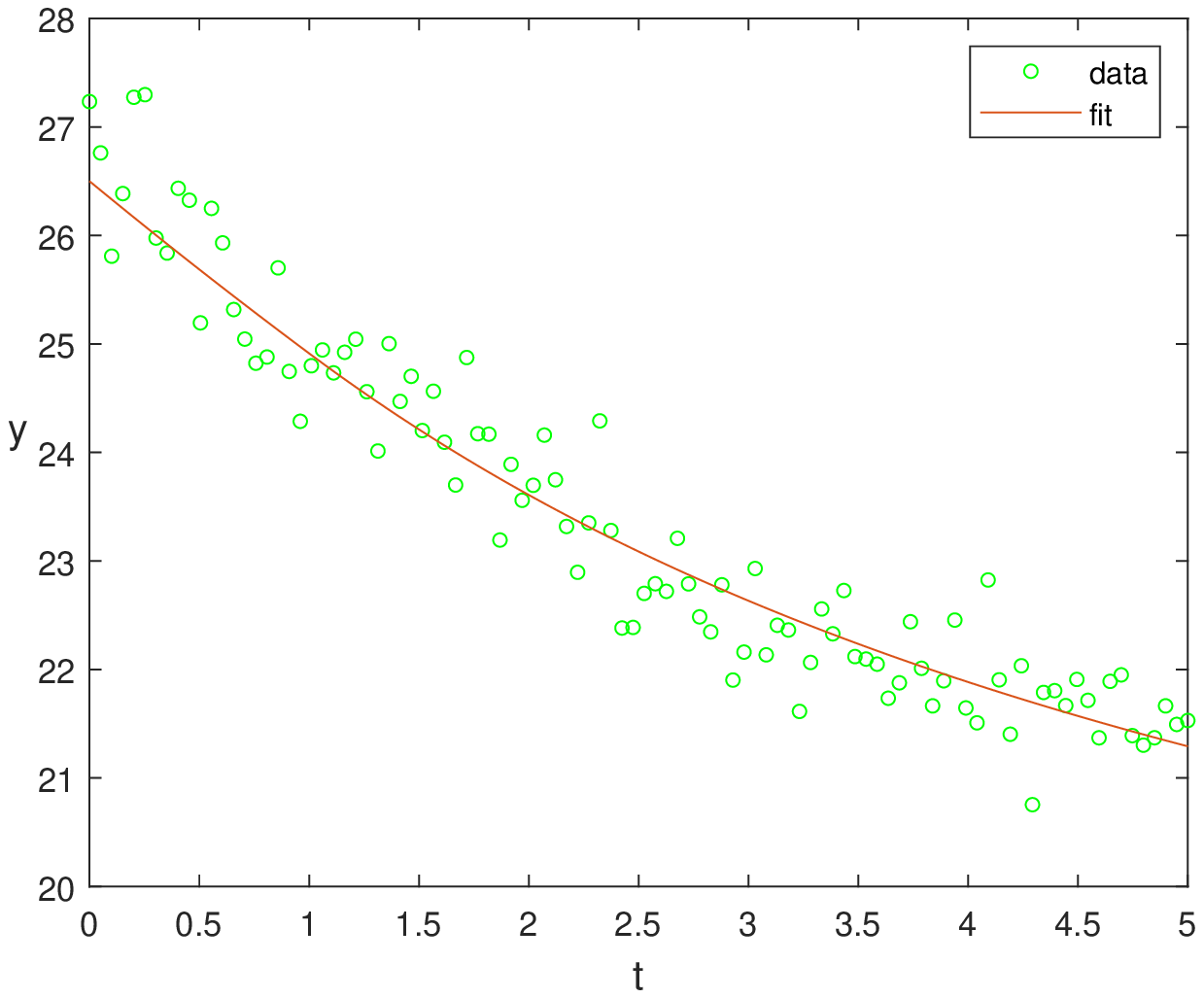}		
		\label{fig:para_est4_try_alpha0p01beta0p01}
	\end{subfigure}
	\caption{Data and curve fit for $x$ and $y$ with Trust-region algorithm for Example 4}
	\label{fig:para_est4_tr_alpha0p01beta0p01}
\end{figure}

\begin{table}[H]
	\begin{center}
		\caption{Estimation of $\alpha$ and $\beta$ with  Trust-region algorithm for Example 4}
		\label{table:ex4_tr}
		\begin{tabular}{ccccr}
			\hline
			Iteration &  \multicolumn{1}{p{2cm}}{\centering Function\\ count} & Residual & Norm of step &  \multicolumn{1}{p{2cm}}{\centering First-order\\ optimality}      \\
			\hline
			0 & 3  & 41221.9 & 1.61e+04    &          \\
			1 & 6  & 9605.73 & 1.95602     & 6.78e+03 \\
			2 & 9  & 1570.02 & 0.941485    & 2.05e+03 \\
			3 & 12 & 171.046 & 0.4507      & 481      \\
			4 & 15 & 40.837  & 0.166794    & 68.6     \\
			5 & 18 & 37.5656 & 0.0302478   & 2.47     \\
			6 & 21 & 37.5613 & 0.00112156  & 0.00392  \\
			7 & 24 & 37.5613 & 1.12797e-06 & 8.64e-06 \\
			\hline
		\end{tabular}
	\end{center}
\end{table}

After seven iterations it gives the best fit of the parameters.

\begin{table}[H]
	\begin{center}
		\caption{True and best fit of  $\alpha$ and $\beta$ for Example 4}
		\label{table:ex4_true_fit_tr}
		\begin{tabular}{|c|c|c|c|}
			\hline
			&Initial & True & Best fit \\
			\hline
			$\alpha$ & 0.3 & 0.5  & 0.5042  \\ \hline
			$\beta$  & 0.5 & 0.8 & 0.7992 \\
			\hline
		\end{tabular}
	\end{center}
\end{table}

\section*{Example 5}
In this example we look at system when it is at the equilibrium. The initial condition are $x(t) = 29.1842$ and $y(t) = 18.2401$ which corresponds to the equilibrium values. For the initial starting point of the parameters, we use $\alpha = 0.01$ and $\beta = 0.01$ 
with random measurement noise $\mathcal{N}(0, 0.20)$ similar to Example 2.

\begin{figure}[H]
	\centering
	\begin{subfigure}[b]{0.48\textwidth}
		\centering
		\includegraphics[width=1\textwidth]{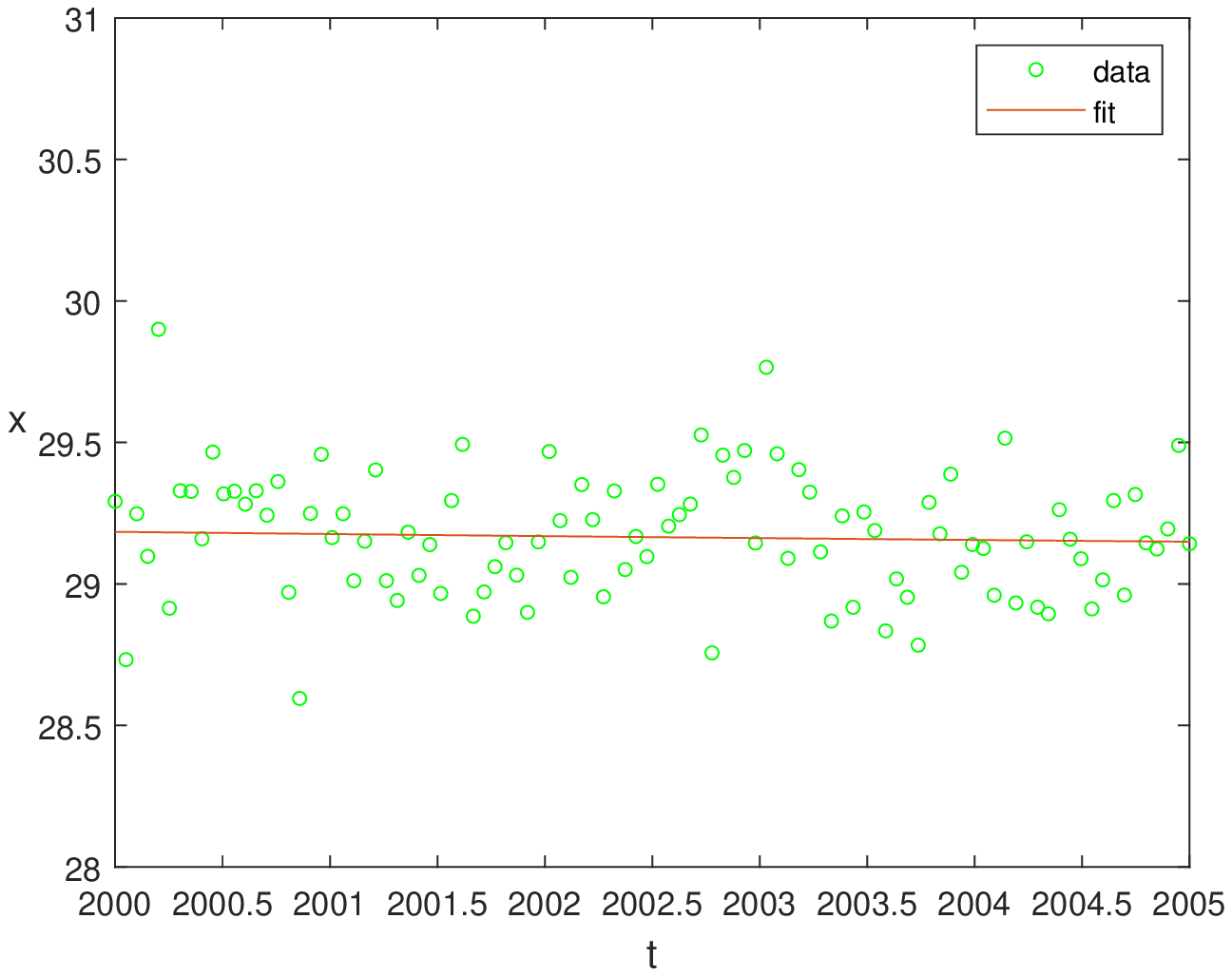}		
		\label{fig:para_est5_lmx_alpha0p01beta0p01}
	\end{subfigure}
	\hfill
	\begin{subfigure}[b]{0.48\textwidth}
		\centering
		\includegraphics[width=1\textwidth]{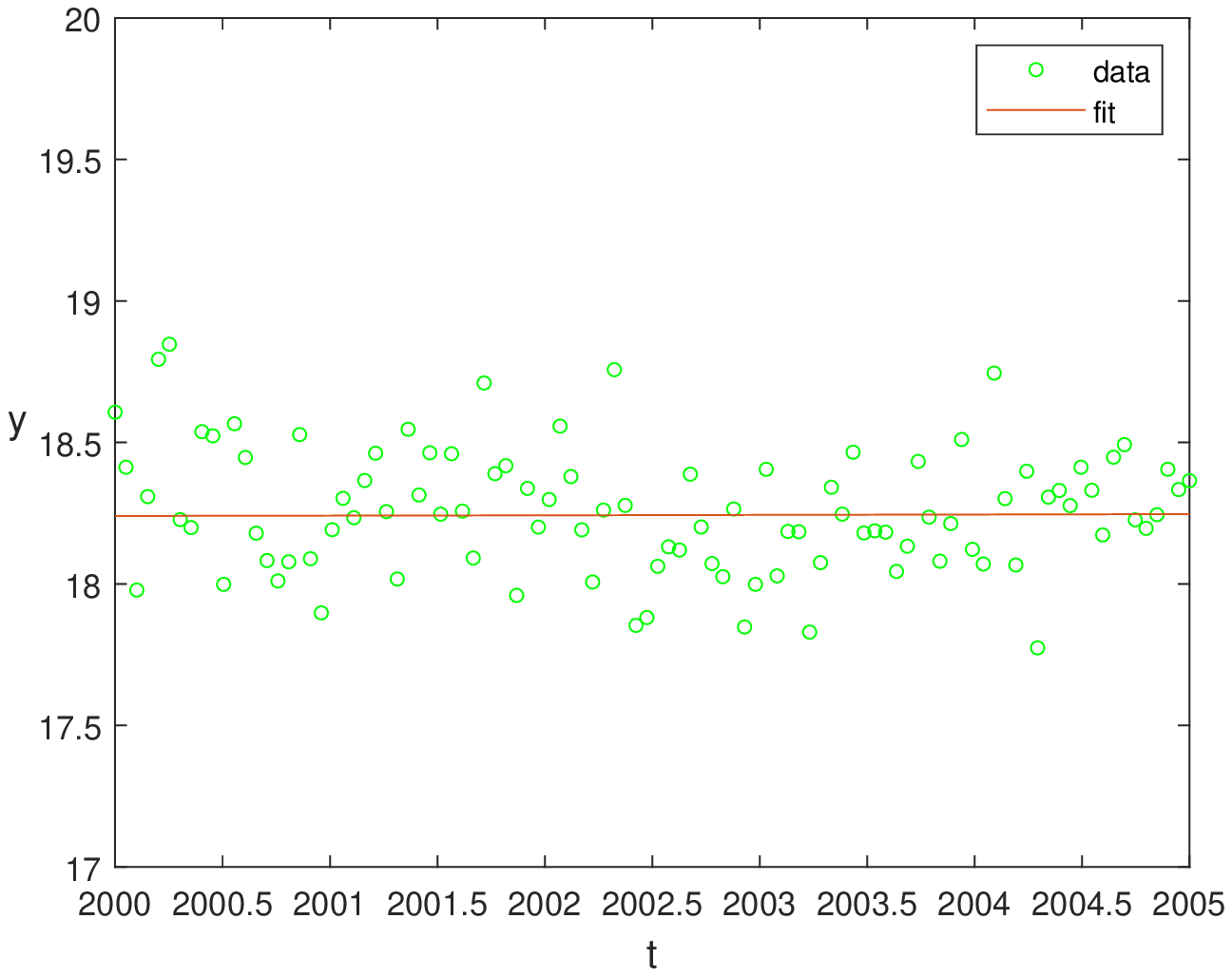}		
		\label{fig:para_est5_lmy_alpha0p01beta0p01}
	\end{subfigure}
	\caption{Data and curve fit for $x$ and $y$ with Levenberg–Marquardt algorithm for Example 5}
	\label{fig:para_est5_lm_alpha0p01beta0p01}
\end{figure}

\begin{table}[H]
	\begin{center}
		\caption{Estimation of $\alpha$ and $\beta$ with  Levenberg–Marquardt algorithm  for Example 5}
		\label{table:ex5_lm}
		\begin{tabular}{cccccr}
			\hline
			Iteration &  \multicolumn{1}{p{2cm}}{\centering Function\\ count} & Residual &  \multicolumn{1}{p{2cm}}{\centering First-order\\ optimality} & Lambda & Norm of step    \\
			\hline
			0 & 3  & 1630.31 & 1.92e+03 & 0.01   &             \\
			1 & 6  & 40.987  & 187      & 0.001  & 0.774995    \\
			2 & 9  & 9.4066  & 3.64     & 0.0001 & 0.152423    \\
			3 & 12 & 9.38722 & 0.00402  & 1e-05  & 0.00420468  \\
			4 & 15 & 9.38722 & 1.91e-06 & 1e-06  & 4.20596e-06 \\
			\hline
		\end{tabular}
	\end{center}
\end{table}

After four iterations it gives the same best fit of the parameters.

\begin{table}[H]
	\begin{center}
		\caption{True and best fit of  $\alpha$ and $\beta$  for Example 5}
		\label{table:ex5_true_fit_lm}
		\begin{tabular}{|c|c|c|c|}
			\hline
			&Initial & True & Best fit \\
			\hline
			$\alpha$ & 0.01 & 0.5  & 0.5040  \\ \hline
			$\beta$  & 0.01 & 0.8 & 0.7990 \\
			\hline
		\end{tabular}
	\end{center}
\end{table}

\begin{figure}[H]
	\centering
	\begin{subfigure}[b]{0.48\textwidth}
		\centering
		\includegraphics[width=1\textwidth]{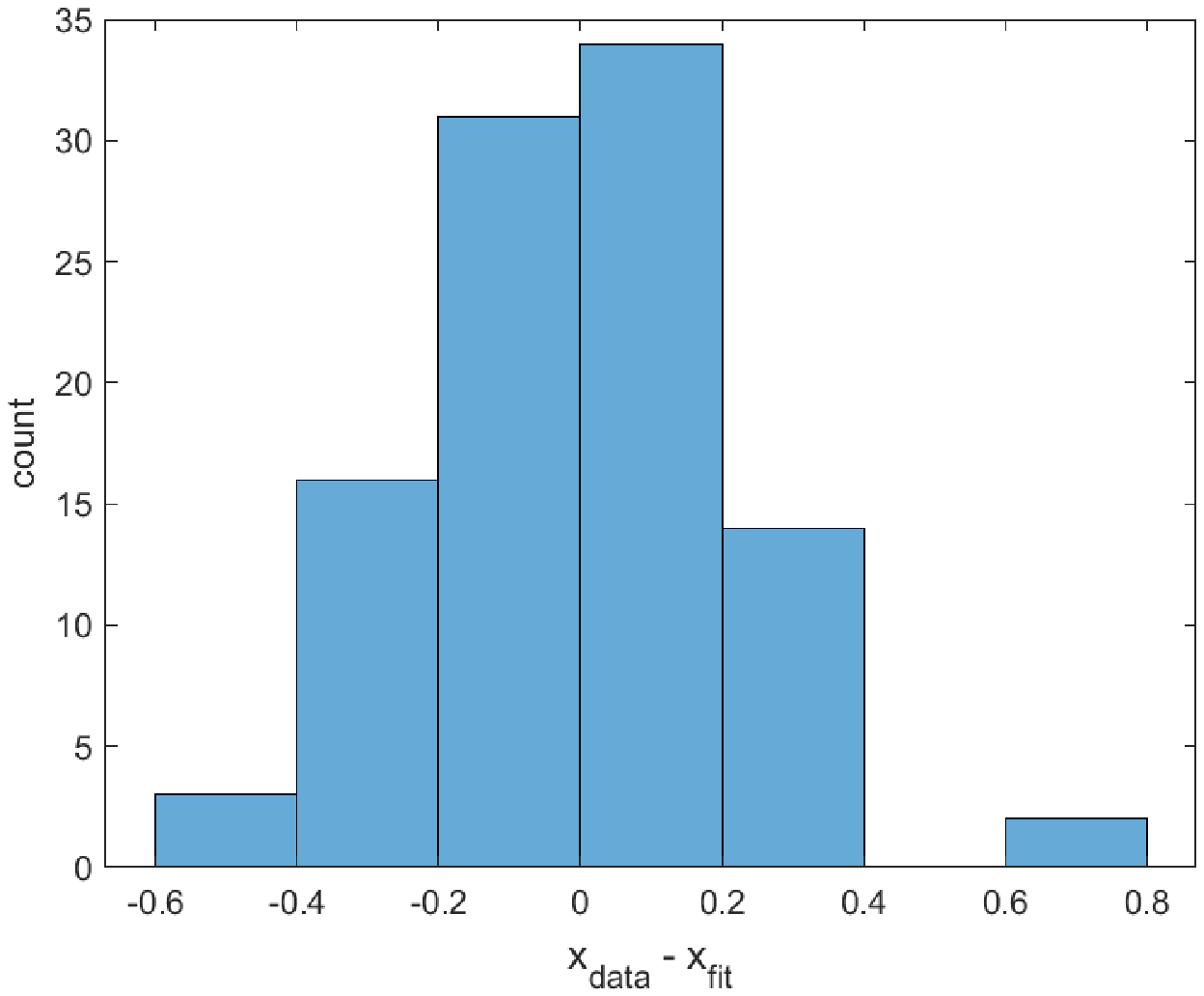}		
		\label{fig:para_est5_hist_lmx_alpha0p01beta0p01}
	\end{subfigure}
	\hfill
	\begin{subfigure}[b]{0.48\textwidth}
		\centering
		\includegraphics[width=1\textwidth]{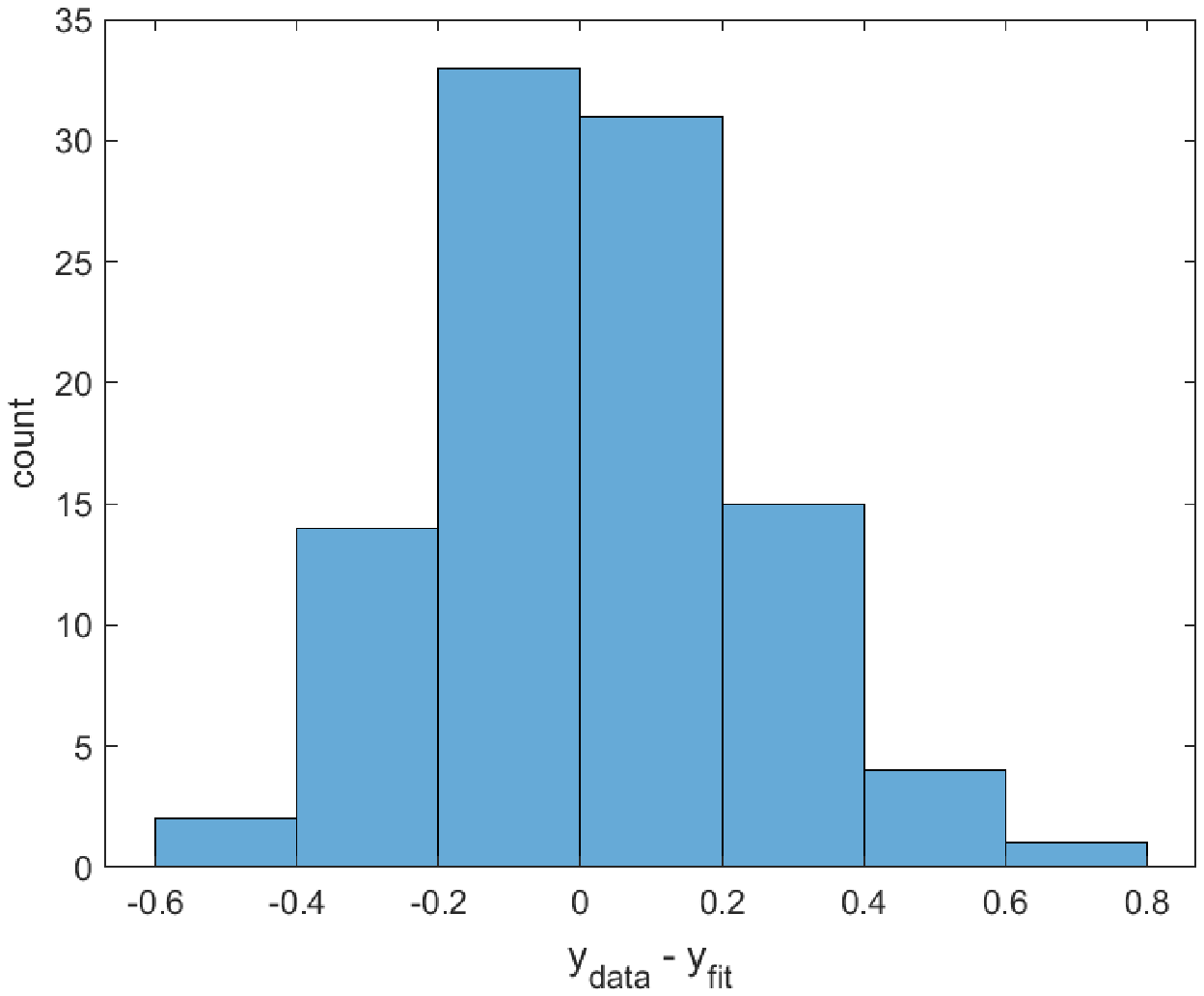}		
		\label{fig:para_est5_hist_lmy_alpha0p01beta0p01}
	\end{subfigure}
	\caption{Histogram of the errors between the data and curve fit with Levenberg–Marquardt algorithm for Example 5}
	\label{fig:para_est5_hist_lm_alpha0p01beta0p01}
\end{figure}

\subsection*{Computation of Estimates with Trust-region algorithm}

\begin{figure}[H]
	\centering
	\begin{subfigure}[b]{0.48\textwidth}
		\centering
		\includegraphics[width=1\textwidth]{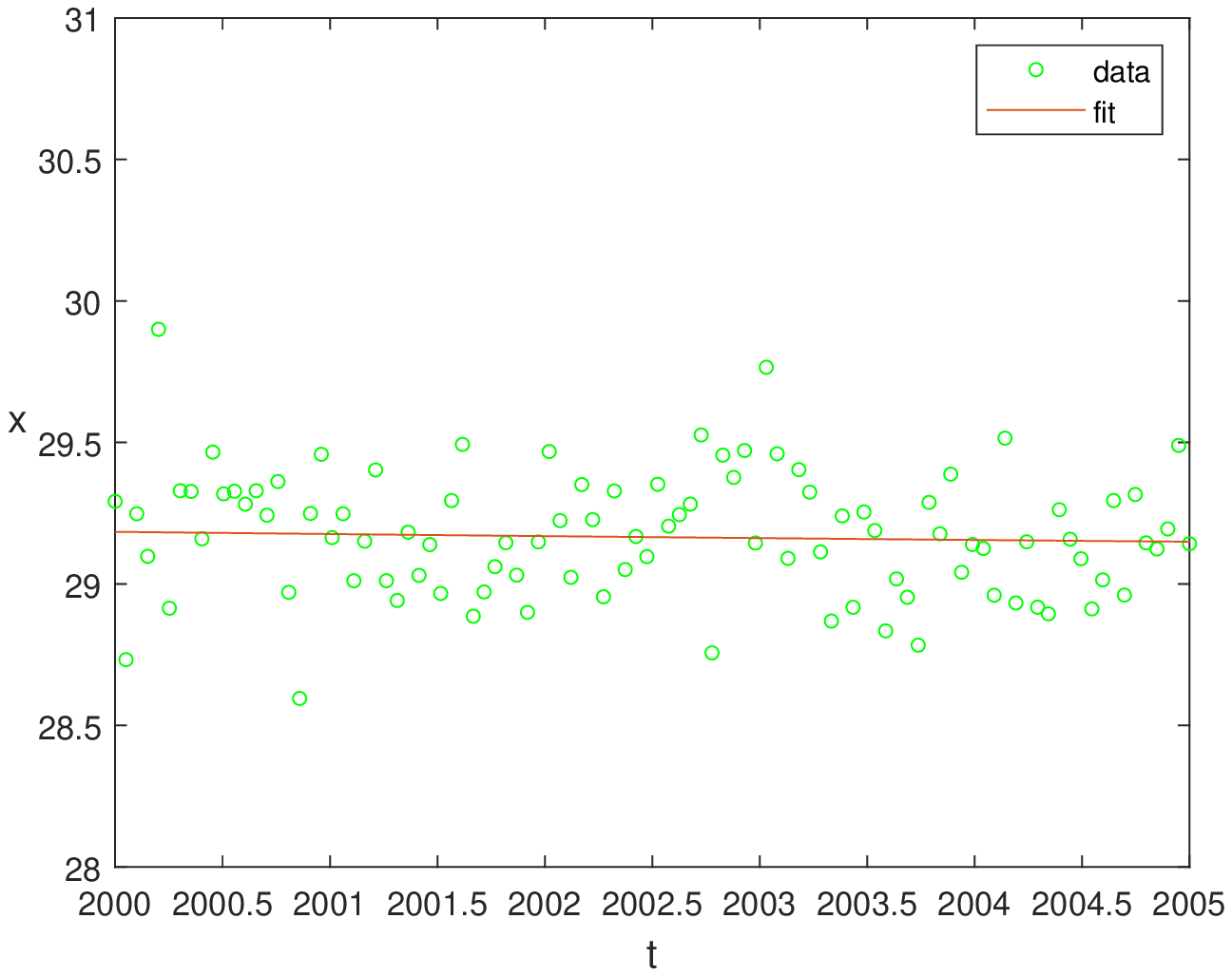}		
		\label{fig:para_est5_trx_alpha0p01beta0p01}
	\end{subfigure}
	\hfill
	\begin{subfigure}[b]{0.48\textwidth}
		\centering
		\includegraphics[width=1\textwidth]{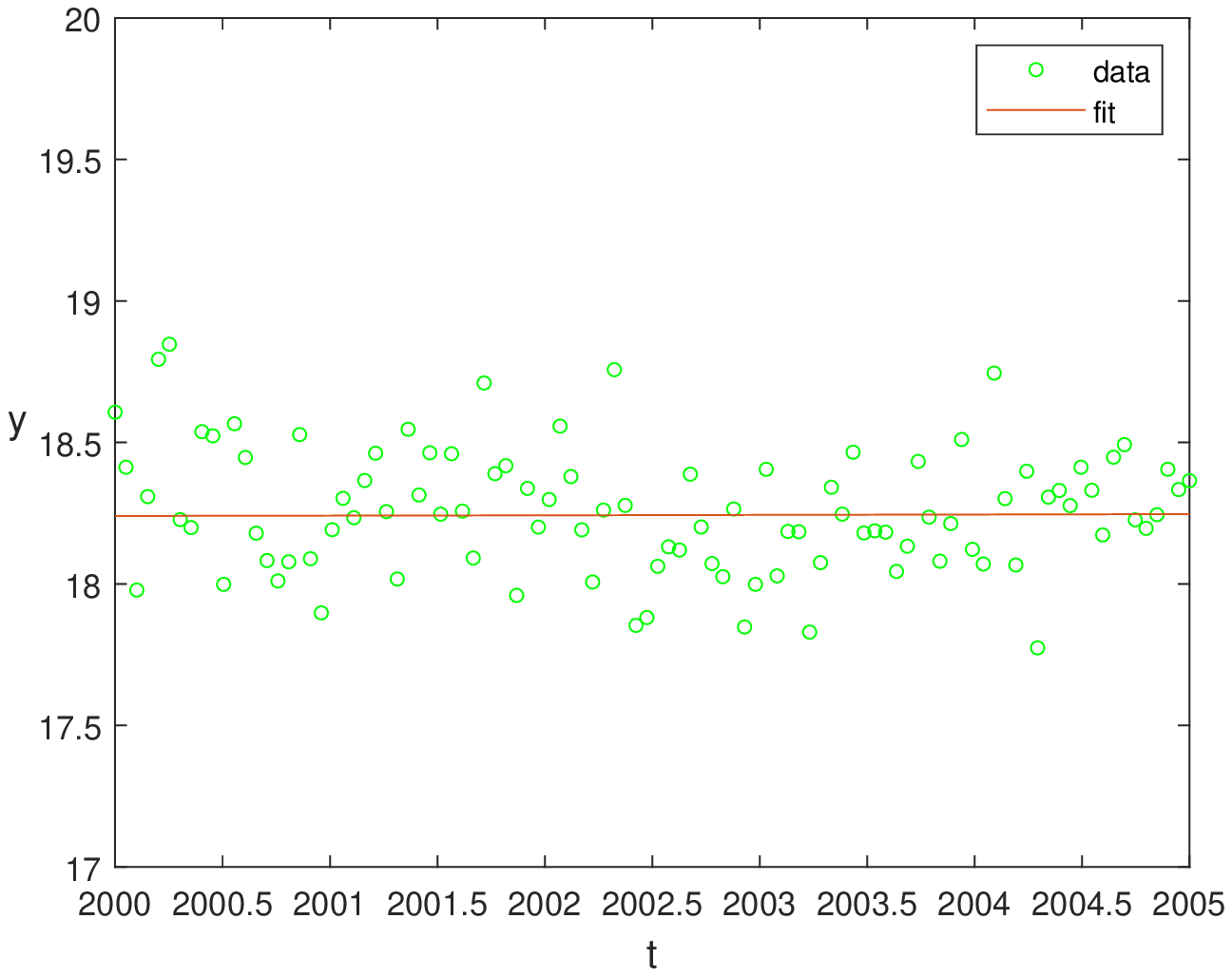}		
		\label{fig:para_est5_try_alpha0p01beta0p01}
	\end{subfigure}
	\caption{Data and curve fit for $x$ and $y$ with Trust-region algorithm for Example 5}
	\label{fig:para_est5_tr_alpha0p01beta0p01}
\end{figure}

\begin{table}[H]
	\begin{center}
		\caption{Estimation of $\alpha$ and $\beta$ with  Trust-region algorithm  for Example 5}
		\label{table:ex5_tr}
		\begin{tabular}{ccccr}
			\hline
			Iteration &  \multicolumn{1}{p{2cm}}{\centering Function\\ count} & Residual & Norm of step &  \multicolumn{1}{p{2cm}}{\centering First-order\\ optimality}      \\
			\hline
			0 & 3  & 18743.5 & 8.96e+03    & 8.96e+03   \\
			1 & 6  & 4288.67 & 1.87713     & 3.32e+03 \\
			2 & 9  & 742.288 & 0.956862    & 954      \\
			3 & 12 & 88.1657 & 0.492103    & 230      \\
			4 & 15 & 12.5793 & 0.208053    & 39.3     \\
			5 & 18 & 9.40059 & 0.0506639   & 2.44     \\
			6 & 21 & 9.38722 & 0.00350099  & 0.0116   \\
			7 & 24 & 9.38722 & 1.66463e-05 & 5.13e-05 \\
			\hline
		\end{tabular}
	\end{center}
\end{table}

After seven iterations it gives the best fit of the parameters.

\begin{table}[H]
	\begin{center}
		\caption{True and best fit of  $\alpha$ and $\beta$  for Example 5}
		\label{table:ex5_true_fit_tr}
		\begin{tabular}{|c|c|c|c|}
			\hline
			&Initial & True & Best fit \\
			\hline
			$\alpha$ & 0.3 & 0.5  & 0.5040  \\ \hline
			$\beta$  & 0.5 & 0.8 & 0.7990 \\
			\hline
		\end{tabular}
	\end{center}
\end{table}

\begin{figure}[H]
	\centering
	\begin{subfigure}[b]{0.48\textwidth}
		\centering
		\includegraphics[width=1\textwidth]{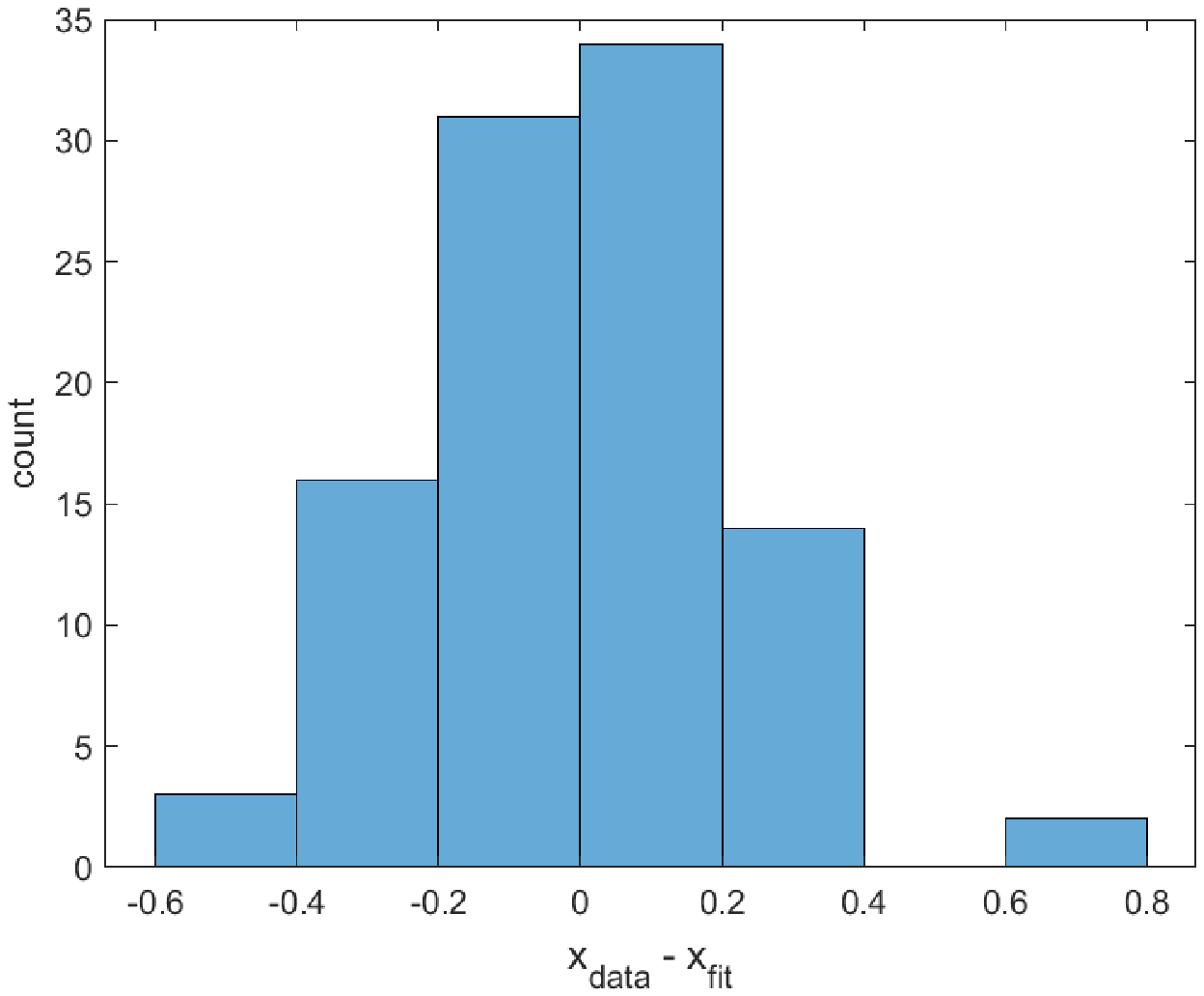}		
		\label{fig:para_est5_hist_trx_alpha0p01beta0p01}
	\end{subfigure}
	\hfill
	\begin{subfigure}[b]{0.48\textwidth}
		\centering
		\includegraphics[width=1\textwidth]{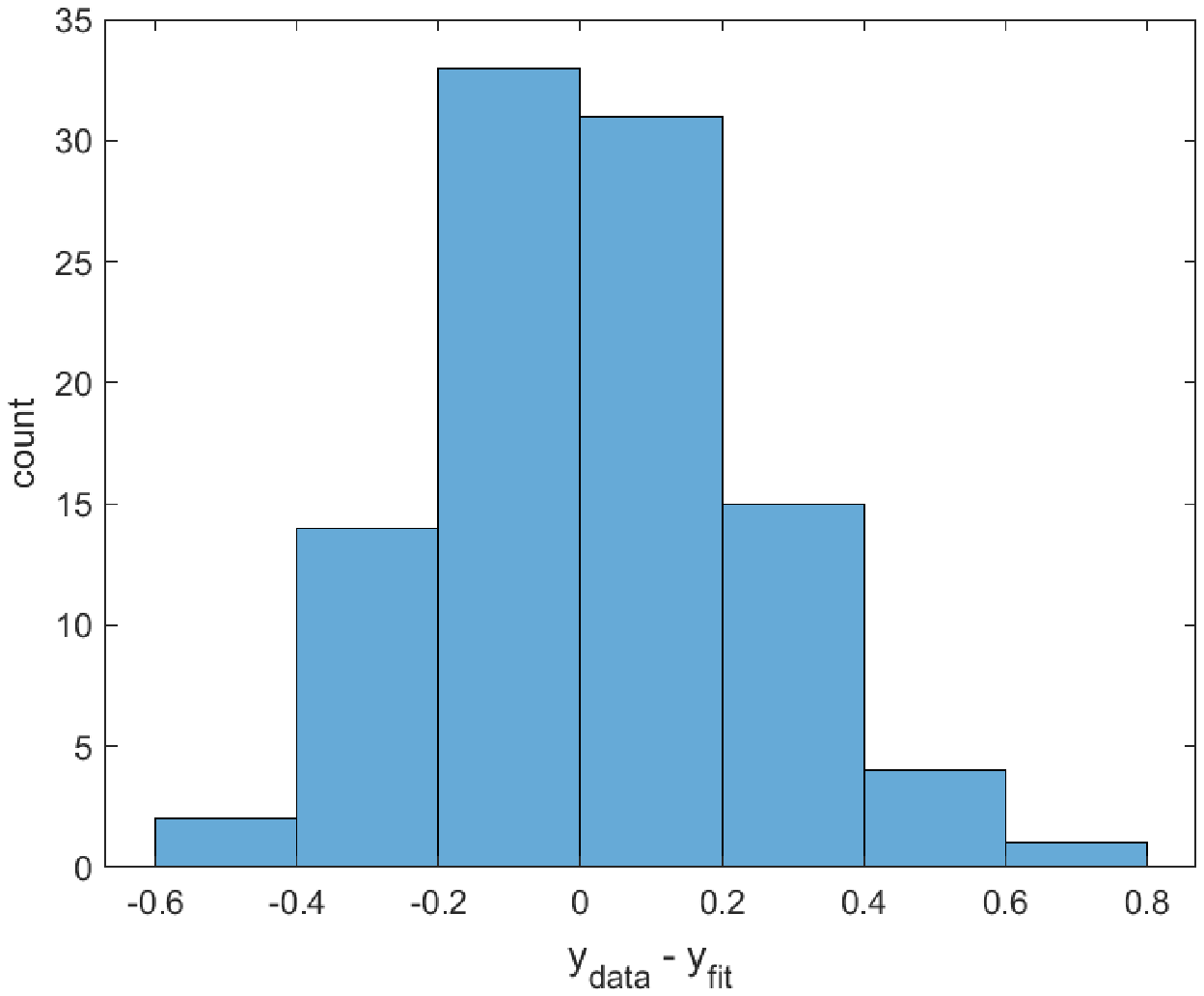}		
		\label{fig:para_est5_hist_try_alpha0p01beta0p01}
	\end{subfigure}
	\caption{Histogram of the errors between the data and curve fit with Trust-region algorithm for Example 5}
	\label{fig:para_est5_hist_tr_alpha0p01beta0p01}
\end{figure}

\section{Summary of the Examples}

\begin{table}[H]
	\begin{center}
		\caption{Summary of  parameter estimation examples where the true parameters are $\alpha = 0.5$ and $\beta = 0.8$}
		\label{table:ex_summary}
		\begin{tabular}{ccccc}
			\hline
			& Initial & Noise &  LM fit &  TR fit      \\
			\hline
			Example 1 & \multicolumn{1}{p{2cm}}{\centering $\alpha = 0.3$\\ $\beta = 0.5$}  & $\mathcal{N}(0, 0.20)$ & \multicolumn{1}{p{2cm}}{\centering $\alpha = 0.5021$\\ $\beta = 0.7996$}   &    \multicolumn{1}{p{2cm}}{\centering $\alpha = 0.5021$\\ $\beta = 0.7996$}      \\ \hline
			Example 2 & \multicolumn{1}{p{2cm}}{\centering $\alpha = 0.01$\\ $\beta = 0.01$}  & $\mathcal{N}(0, 0.20)$ & \multicolumn{1}{p{2cm}}{\centering $\alpha = 0.5021$\\ $\beta = 0.7996$}   &    \multicolumn{1}{p{2cm}}{\centering $\alpha = 0.5021$\\ $\beta = 0.7996$}      \\ \hline
			Example 3 & \multicolumn{1}{p{2cm}}{\centering $\alpha = 0.3$\\ $\beta = 0.5$}  & $\mathcal{N}(0, 0.40)$ & \multicolumn{1}{p{2cm}}{\centering $\alpha = 0.5042$\\ $\beta = 0.7992$}   &    \multicolumn{1}{p{2cm}}{\centering $\alpha = 0.5042$\\ $\beta = 0.7992$}      \\ \hline
			Example 4 & \multicolumn{1}{p{2cm}}{\centering $\alpha = 0.01$\\ $\beta = 0.01$}  & $\mathcal{N}(0, 0.40)$ & \multicolumn{1}{p{2cm}}{\centering $\alpha = 0.5042$\\ $\beta = 0.7992$}   &    \multicolumn{1}{p{2cm}}{\centering $\alpha = 0.5042$\\ $\beta = 0.7992$}      \\ \hline
			Example 5 & \multicolumn{1}{p{2cm}}{\centering $\alpha = 0.01$\\ $\beta = 0.01$}  & $\mathcal{N}(0, 0.20)$ & \multicolumn{1}{p{2cm}}{\centering $\alpha = 0.5040$\\ $\beta = 0.7990$}   &    \multicolumn{1}{p{2cm}}{\centering $\alpha = 0.5040$\\ $\beta = 0.7990$}      \\
			\hline
		\end{tabular}
	\end{center}
\end{table}

We observe that the two algorithms produce similar results in these simulated examples. The initial starting point and noise didn't have significant impact in the number of iteration it takes to compute the estimates. We were able to compute the estimates both at the starting time interval (Example 1-4) and at the equilibrium (Example 5).

	\bibliography{references1}
\end{document}